\documentclass{amsart}
\usepackage{amsfonts}
\usepackage{amscd}

\input{tcilatex}
\theoremstyle{definition}
\theoremstyle{remark}
\numberwithin{equation}{section}

\input tcilatex

\begin{document}
\title{Infinite Dimensional It\^{o} Algebras\\
of Quantum White Noise}
\author{V. P. Belavkin.}
\address{Mathematics Department, University of Nottingham,\\
NG7 2RD, UK.}
\email{vpb@@maths.nott.ac.uk}
\date{January 20, 1998. Extended from: V. P. Belavkin, On Quantum It\^{o}
Algebras and Their Decompositions. Letters in Math Phys \textbf{7}, 1-16,
1998.}
\subjclass{Quantum Probability and Stochastic Analysis}
\thanks{Published in: \textit{Trends in Contemporary Infinite Dimensional
Analysis and Quantum Probability} 57--80, Instituto Italiano di Cultura,
Kyoto, 2000.}

\begin{abstract}
A simple axiomatic characterization of the general (infinite dimensional,
noncommutative) It\^{o} algebra is given and a pseudo-Euclidean fundamental
representation for such algebra is described. The notion of It\^{o}
B*-algebra, generalizing the C*-algebra is defined to include the Banach
infinite dimensional It\^{o} algebras of quantum Brownian and quantum L\'{e}%
vy motion, and the B*-algebras of vacuum and thermal quantum noise are
characterized. It is proved that every It\^{o} algebra is canonically
decomposed into the orthogonal sum of quantum Brownian (Wiener) algebra and
quantum L\'{e}vy (Poisson) algebra. In particular, every quantum thermal
noise is the orthogonal sum of a quantum Wiener noise and a quantum Poisson
noise as it is stated by the L\'{e}vy-Khinchin theorem in the classical case.
\end{abstract}

\maketitle

\section{Introduction}

The classical differential calculus for the infinitesimal increments $%
\mathrm{d}x=x\left( t+\mathrm{d}t\right) -x\left( t\right) $ became
generally accepted only after Newton gave a very simple algebraic rule $%
\left( \mathrm{d}t\right) ^{2}=0$ for the formal computations of first order
differentials for smooth trajectories $t\mapsto x\left( t\right) $ in a
phase space. The linear space of the differentials $\mathrm{d}x=\alpha 
\mathrm{d}t$ for a (complex) trajectory became treated at each $x=x\left(
t\right) \in \mathbb{C}$ as a one-dimensional algebra $\mathfrak{a}=\mathbb{C%
}d_{t}$ of the elements $a=\alpha d_{t}$ with involution $a^{\star }=\bar{%
\alpha}d_{t}$ given by the complex conjugation $\alpha \mapsto \bar{\alpha}$
of the derivative $\alpha =\mathrm{d}x/\mathrm{d}t\in \mathbb{C}$ and the
nilpotent multiplication $a\cdot a^{\star }=0$ corresponding to the
multiplication table $d_{t}\cdot d_{t}^{\star }=0$ for the basic nilpotent
element $d_{t}=d_{t}^{\star }$, the abstract notation of $\mathrm{d}t$. Note
that the nilpotent $\star $-algebra $\mathfrak{a}$ of abstract
infinitesimals $\alpha d_{t}$ has no realization in complex numbers, as well
as no operator representation $\alpha D_{t}$ with a Hermitian nilpotent $%
D_{t}=D_{t}^{\dagger }$ in a Euclidean (complex pre-Hilbert) space, but it
can be represented by the algebra of complex nilpotent $2\times 2$ matrices $%
\hat{a}=\alpha \hat{d}_{t}$, where $\hat{d}_{t}=\frac{1}{2}\left( \hat{\sigma%
}_{3}+i\hat{\sigma}_{1}\right) =\hat{d}_{t}^{\dagger }$ with respect to the
standard Minkowski metric $\left( \mathrm{x}|\mathrm{x}\right) =\left\vert
\zeta \right\vert ^{2}-\left\vert \eta \right\vert ^{2}$ for $\mathrm{x}%
=\zeta \mathrm{e}_{+}+\eta \mathrm{e}_{-}$ in $\mathbb{C}^{2}$. The complex
pseudo-Hermitian nilpotent matrix $\hat{d}_{t}$, $\hat{d}_{t}^{2}=0$,
representing the multiplication $d_{t}^{2}=d_{t}\cdot d_{t}=0$, has the
canonical triangular form 
\begin{equation}
\text{$\mathrm{D}_{t}$}=\left[ 
\begin{array}{ll}
0 & 1 \\ 
0 & 0%
\end{array}%
\right] ,\quad \quad \text{$\mathrm{D}_{t}$}^{\ast }=\left[ 
\begin{array}{ll}
0 & 0 \\ 
1 & 0%
\end{array}%
\right] ,\quad
\end{equation}%
\begin{equation*}
\text{$\mathrm{D}_{t}$}^{\dagger }=\left[ 
\begin{array}{ll}
0 & 1 \\ 
1 & 0%
\end{array}%
\right] \text{$\mathrm{D}_{t}$}^{\ast }\left[ 
\begin{array}{ll}
0 & 1 \\ 
1 & 0%
\end{array}%
\right] =\text{$\mathrm{D}_{t}$}
\end{equation*}%
in the basis $\mathrm{h}_{\pm }=\left( \mathrm{e}_{+}\pm \mathrm{e}%
_{-}\right) /\sqrt{2}$ in which $\left( \mathbf{x}\text{$\mathrm{D}_{t}$}|%
\mathbf{x}\right) =\overline{\left( \mathbf{x}\text{$\mathrm{D}_{t}$}|%
\mathbf{x}\right) }$ for all $\mathbf{x}=\left( \xi _{-},\xi _{+}\right) $
with respect to the pseudo-Euclidean scalar product $\left( \mathbf{x}|%
\mathbf{x}\right) =\xi _{-}\xi ^{-}+\xi _{+}\xi ^{+}$, where $\xi ^{\pm
}=\left( \zeta \pm \eta \right) /\sqrt{2}=\bar{\xi}_{\mp }\in \mathbb{C}$.

The Newton's formal computations can be generalized to non-smooth paths to
include the calculus of first order forward differentials $\mathrm{d}y\simeq
\left( \mathrm{d}t\right) ^{1/2}$ of continuous diffusions $y\left( t\right)
\in \mathbb{R}$ which have no derivative at any $t$, and the forward
differentials $\mathrm{d}n\in \left\{ 0,1\right\} $ of left continuous
counting trajectories $n\left( t\right) \in \mathbb{Z}_{+}$ which have zero
derivative for almost all $t$ (except the points of discontinuity when $%
\mathrm{d}n=1$). The first is usually done by adding the rules 
\begin{equation}
\left( \mathrm{d}w\right) ^2=\mathrm{d}t,\quad \mathrm{d}w\mathrm{d}t=0=%
\mathrm{d}t\mathrm{d}w  \label{2}
\end{equation}
in formal computations of continuous trajectories having the first order
forward differentials $\mathrm{d}x=\alpha \mathrm{d}t+\beta \mathrm{d}w$
with the diffusive part given by the increments of standard Brownian paths $%
w\left( t\right) $. The second can be done by adding the rules 
\begin{equation}
\left( \mathrm{d}m\right) ^2=\mathrm{d}m+\mathrm{d}t,\quad \mathrm{d}m%
\mathrm{d}t=0=\mathrm{d}t\mathrm{d}m  \label{3}
\end{equation}
in formal computations of left continuous and smooth for almost all $t$
trajectories having the forward differentials $\mathrm{d}x=\alpha \mathrm{d}%
t+\gamma \mathrm{d}m$ with jumping part $\mathrm{d}z\in \left\{ \gamma
,-\gamma \mathrm{d}t\right\} $ given by the increments of standard L\'{e}vy
paths $m\left( t\right) =n\left( t\right) -t$. These rules, well known since
the beginning of this century, were formalized by It\^{o} \cite{Ito51} into
the form of a stochastic calculus: the first one is now known as the
multiplication rule for the forward differential of the standard Wiener
process $w\left( t\right) $, and the second one is the multiplication rule
for the forward differential of the standard Poisson process $n\left(
t\right) $, compensated by its mean value $t$.

The linear span of $\mathrm{d}t$ and $\mathrm{d}w$ forms a two-dimensional
differential It\^{o} algebra $\mathfrak{b}=\mathbb{C}d_{t}+\mathbb{C}d_{w}$
for the complex Brownian motions $x\left( t\right) =\int \alpha \mathrm{d}%
t+\int \eta \mathrm{d}w$, where $d_{w}=d_{w}^{\star }$ is a nilpotent of
second order element, representing the real increment $\mathrm{d}w$, with
multiplication table $d_{w}^{2}=d_{t}$, $d_{w}\cdot d_{t}=0=d_{t}\cdot d_{w}$%
, while the linear span of $\mathrm{d}t$ and $\mathrm{d}m$ forms a
two-dimensional differential It\^{o} algebra $\mathfrak{c}=\mathbb{C}d_{t}+%
\mathbb{C}d_{m}$ for the complex L\'{e}vy motions $x=\int \alpha \mathrm{d}%
t+\int \zeta \mathrm{d}m$, where $d_{m}=d_{m}^{\star }$ is a basic element,
representing the real increment $\mathrm{d}m$, with multiplication table $%
d_{m}^{2}=d_{m}+d_{t}$, $d_{m}\cdot d_{t}=0=d_{t}\cdot d_{m}$. As in the
case of the Newton algebra, the It\^{o} $\star $-algebras $\mathfrak{b}$ and 
$\mathfrak{c}$ have no Euclidean operator realization, but they can be
represented by the algebras of triangular matrices $\mathrm{B}=\alpha 
\mathrm{D}_{t}+\eta \mathrm{D}_{w}$, $\mathrm{C}=\alpha \mathrm{D}_{t}+\zeta 
\mathrm{D}_{m}$ with pseudo-Hermitian basis elements 
\begin{equation}
\text{$\mathrm{D}_{t}$}=\left[ 
\begin{array}{lll}
0 & 0 & 1 \\ 
0 & 0 & 0 \\ 
0 & 0 & 0%
\end{array}%
\right] =\text{$\mathrm{D}_{t}$}^{\dagger },\quad \mathrm{D}_{w}=\left[ 
\begin{array}{lll}
0 & 1 & 0 \\ 
0 & 0 & 1 \\ 
0 & 0 & 0%
\end{array}%
\right] =\mathrm{D}_{w}^{\dagger },  \label{4}
\end{equation}%
\begin{equation*}
\emph{\quad }\mathrm{D}_{m}\mathbf{=}\left[ 
\begin{array}{lll}
0 & 1 & 0 \\ 
0 & 1 & 1 \\ 
0 & 0 & 0%
\end{array}%
\right] =\mathrm{D}_{m}^{\dagger },
\end{equation*}%
where $\left( \mathbf{x}\mathrm{B}^{\dagger }|\mathbf{x}\right) =\overline{%
\left( \mathbf{x}\mathrm{B}|\mathbf{x}\right) }$ for all $\mathbf{x}=\left(
\xi _{-},\xi _{\circ },\xi _{+}\right) \in \mathbb{C}^{3}$ in the complex
three-dimensional Minkowski space with respect to the indefinite metric $%
\left( \mathbf{x}|\mathbf{x}\right) =\xi _{-}\xi ^{-}+\xi _{\circ }\xi
^{\circ }+\xi _{+}\xi ^{+}$, where $\xi ^{\mu }=\bar{\xi}_{-\mu }$ with $%
-\left( -,\circ ,+\right) =\left( +,\circ ,-\right) $.

Note that according to the L\'{e}vy-Khinchin theorem, every stochastic
process $x\left( t\right) $ with independent increments can be canonically
decomposed into a smooth, Wiener and Poisson parts as in the mixed case of
one-dimensional complex motion $x\left( t\right) =\int \alpha \mathrm{d}%
t+\int \eta \mathrm{d}w+\int \zeta \mathrm{d}m$ given by the orthogonal and
thus commutative increments $\mathrm{d}w\mathrm{d}m=0=\mathrm{d}m\mathrm{d}w$%
. In fact such generalized commutative differential calculus applies not
only to the stochastic integration with respect to the processes with
independent increments; these formal algebraic rules, or their
multidimensional versions, can be used for formal computations of forward
differentials for any classical trajectories decomposed into the smooth,
diffusive and jumping parts.

Two natural questions arrise: are there other then these two commutative
differential algebras which could be useful, in particular, for formal
computations of the noncommutative differentials in quantum theory, and if
there are, is it possible to characterize them by simple axioms and to give
a generalized version of the L\'{e}vy-Khinchin decomposition theorem? The
first question has been already positively answered since the well known
differential realization of the simplest non-commutative table $%
d_{w}d_{w}^{\star }=\rho _{+}d_{t}$, $d_{w}^{\star }d_{w}=\rho _{-}d_{t}$\
for $\rho _{+}>\rho _{-}\geq 0$ was given in the mid of 60-th in terms of
the annihilators $\hat{w}\left( t\right) $ and creators $\hat{w}^{\dagger
}\left( t\right) $ of a quantum Brownian thermal noise \cite{Grd91}. This
paper gives a systematic answer on the second question, the first part of
which has been in principle positively resolved in our papers \cite{3,4}.

Although the orthogonality condition $d_{w}\cdot d_{m}=0=d_{w}\cdot d_{m}$
for the classical independent increments $\mathrm{d}w$ and $\mathrm{d}m$ can
be realized only in a higher, at least four, dimensional Minkowski space, it
is interesting to make sense of the non-commutative $\star $-algebra,
generated by three dimensional non-orthogonal matrix representations (\ref{4}%
) of these differentials with $d_{w}\cdot d_{m}\neq d_{w}\cdot d_{m}$: 
\begin{eqnarray*}
\mathrm{D}_{w}\mathrm{D}_{m} &=&\left( \mathrm{D}_{m}\mathrm{D}_{w}\right)
^{\dagger }=\left[ 
\begin{array}{lll}
0 & 1 & 1 \\ 
0 & 0 & 0 \\ 
0 & 0 & 0%
\end{array}%
\right] \\
&\neq &\left[ 
\begin{array}{lll}
0 & 0 & 1 \\ 
0 & 0 & 1 \\ 
0 & 0 & 0%
\end{array}%
\right] =\left( \mathrm{D}_{w}\mathrm{D}_{m}\right) ^{\dagger }=\mathrm{D}%
_{m}\mathrm{D}_{w}.
\end{eqnarray*}%
This is the four-dimensional $\star $-algebra $\mathfrak{a}=\mathbb{C}%
\mathrm{D}_{t}+\mathbb{C}\mathrm{E}_{-}+\mathbb{C}\mathrm{E}^{+}+\mathbb{C}%
\mathrm{E}$ of triangular matrices $\mathrm{A}=\alpha \mathrm{D}+z^{-}%
\mathrm{E}_{-}+z_{+}\mathrm{E}^{+}+z\mathrm{E}$, where $\mathrm{E}^{+}=%
\mathrm{GE}_{-}^{\ast }\mathrm{G}=\mathrm{E}_{-}^{\dagger }$, $\mathrm{E}=%
\mathrm{E}^{\dagger }$ with respect to the Minkowski metric tensor $\mathrm{G%
}$ in the canonical basis, 
\begin{equation*}
\mathrm{G}=\left[ 
\begin{array}{lll}
0 & 0 & 1 \\ 
0 & 1 & 0 \\ 
1 & 0 & 0%
\end{array}%
\right] ,\,\mathrm{E}_{-}=\left[ 
\begin{array}{lll}
0 & 1 & 0 \\ 
0 & 0 & 0 \\ 
0 & 0 & 0%
\end{array}%
\right] ,\,\mathrm{E}^{+}\mathbf{=}\left[ 
\begin{array}{lll}
0 & 0 & 0 \\ 
0 & 0 & 1 \\ 
0 & 0 & 0%
\end{array}%
\right] ,\,\mathrm{E}=\left[ 
\begin{array}{lll}
0 & 0 & 0 \\ 
0 & 1 & 0 \\ 
0 & 0 & 0%
\end{array}%
\right] ,
\end{equation*}%
given by the algebraic combinations 
\begin{equation*}
\mathrm{E}_{-}=\mathrm{D}_{w}\mathrm{D}_{m}-\mathrm{D}_{t},\mathrm{D}^{+}=%
\mathrm{D}_{m}\mathrm{D}_{w}-\mathrm{D}_{t},\mathrm{E}=\mathrm{D}_{m}-%
\mathrm{D}_{w}
\end{equation*}%
of three matrices (\ref{4}). It realizes the multiplication table 
\begin{equation*}
e_{-}\cdot e^{+}=\text{$d_{t}$},\quad e_{-}\cdot e=e_{-},\quad e\cdot
e^{+}=e^{+},\quad e^{2}=e
\end{equation*}%
with the products for all other pairs being zero, unifying the commutative
tables (\ref{2}), (\ref{3}). It is well known in the quantum stochastic
calculus as the HP (Hudson-Parthasarathy) table \cite{1} 
\begin{equation*}
\mathrm{d}\Lambda _{-}\mathrm{d}\Lambda ^{+}=I\mathrm{d}t,\quad \mathrm{d}%
\Lambda _{-}\mathrm{d}\Lambda =\mathrm{d}\Lambda _{-},\quad \mathrm{d}%
\Lambda \mathrm{d}\Lambda ^{+}=\mathrm{d}\Lambda ^{+},\quad \left( \mathrm{d}%
\Lambda \right) ^{2}=\mathrm{d}\Lambda ,
\end{equation*}%
with zero products for all other pairs, for the multiplication of the
canonical number $\mathrm{d}\Lambda $, creation $\mathrm{d}\Lambda ^{+}$,
annihilation $\mathrm{d}\Lambda _{-}$, and preservation $\mathrm{d}\Lambda
_{+}^{-}=I\mathrm{d}t$ differentials in Fock space over the Hilbert space $%
L^{2}\left( \mathbb{R}_{+}\right) $ of square-integrable complex functions $%
f\left( t\right) ,t\in \mathbb{R}_{+}$.

Note that any two-dimensional It\^{o} $\star $-algebra $\mathfrak{a}$ is
commutative as $d_{t}a=0=ad_{t}$ for any other element $a\neq d_{t}$ of the
basis $\left\{ a,d_{t}\right\} $ in $\mathfrak{a}$. Moreover, each such
algebra is either of the Wiener or of the Poisson type, as it is either
second order nilpotent, or contains a unital one-dimensional subalgebra, as
the cases of the subalgebras $\mathfrak{b},\mathfrak{c}$. Other
two-dimensional sub-algebras containing $d_{t}$, are generated by either
Wiener $d_{w}=\bar{\xi}e_{-}+\xi e^{+}$ or Poisson $d_{m}=e+\lambda d_{w}$
element with the special case $d_{m}=e$, corresponding to the only
non-faithful It\^{o} algebra of the Poisson process with zero intensity $%
\lambda ^{2}=0$. However there is only one three dimensional $\star $%
-subalgebra of the four-dimensional HP algebra with $d_{t}$, namely the
noncommutative subalgebra of vacuum Brownian motion, generated by the
creation $e^{+}$ and annihilation $e_{-}$ differentials$.$ Thus our results
on the classification of noncommutative It\^{o} $\star $-algebras will be
nontrivial only in the higher dimensions of $\mathfrak{a}$.

The well known L\'{e}vy-Khinchin classification of the classical noise can
be reformulated in purely algebraic terms as the decomposability of any
commutative It\^{o} algebra into Wiener (Brownian) and Poisson (L\'{e}vy)
orthogonal components. In the general case we shall show that every It\^{o} $%
\star $-algebra is also decomposable into a quantum Brownian, and a quantum L%
\'{e}vy orthogonal components.

Thus classical stochastic calculus developed by It\^{o}, and its quantum
stochastic analog, given by Hudson and Parthasarathy in \cite{1}, has been
unified in a $\star $-algebraic approach to the operator integration in Fock
space \cite{3}, in which the classical and quantum calculi become
represented as two extreme commutative and noncommutative cases of a
generalized It\^{o} calculus.

In the next section we remind the definition of the general It\^{o} algebra
and show that every such algebra can be embedded as a $\star $-subalgebra
into in general infinite dimensional Hudson-Parthasarathy algebra as it was
first proved in \cite{4}.

\section{Representations of It\^{o} $\star $-algebras}

The generalized It\^{o} algebra was defined in \cite{3} as a linear span of
the differentials 
\begin{equation*}
\mathrm{d}\Lambda \left( t,a\right) =\Lambda \left( t+\mathrm{d}t,a\right)
-\Lambda \left( t,a\right) ,\quad a\in \mathfrak{a}
\end{equation*}%
for a family $\left\{ \Lambda \left( a\right) :a\in \mathfrak{a}\right\} $
of operator-valued integrators $\Lambda \left( t,a\right) $ on a pre-Hilbert
space, satisfying for each $t\in \mathbb{R}_{+}$ the $\star $-semigroup
conditions 
\begin{equation}
\Lambda \left( t,a^{\star }\right) =\Lambda \left( t,a\right) ^{\dagger
},\quad \mathrm{d}\Lambda \left( t,a\cdot b\right) =\mathrm{d}\Lambda \left(
t,a\right) \mathrm{d}\Lambda \left( t,b\right) ,\quad \Lambda \left(
t,d_{t}\right) =tI,  \label{0.1}
\end{equation}%
with mean values $\langle \mathrm{d}\Lambda \left( t,a\right) \rangle
=l\left( a\right) \mathrm{d}t$ in a given vector state $\langle \mathrm{%
\cdot }\rangle $, absolutely continuous with respect to $\mathrm{d}t$. Here $%
\Lambda \left( t,a\right) ^{\dagger }$ means the Hermitian conjugation of
the (unbounded) operator $\Lambda \left( t,a\right) $, which is defined on
the pre-Hilbert space for each $t\in \mathbb{R}_{+}$ as the operator $%
\Lambda \left( t,a^{\star }\right) $, 
\begin{eqnarray*}
&&\mathrm{d}\Lambda \left( t,a\right) \mathrm{d}\Lambda \left( t,b\right) \\
&=&\mathrm{d}\left( \Lambda \left( t,a\right) \Lambda \left( t,b\right)
\right) -\mathrm{d}\Lambda \left( t,a\right) \Lambda \left( t,b\right)
-\Lambda \left( t,a\right) \mathrm{d}\Lambda \left( t,b\right) ,
\end{eqnarray*}%
and $\mathrm{d}t$ is embedded into the family of the operator-valued
differentials as $\mathrm{d}\Lambda \left( t,d_{t}\right) $ with the help of
a special element $d_{t}=d_{t}^{\star }$ of the parametrizing $\star $%
-semigroup $\mathfrak{a}$.

Assuming that the parametrization is exact such that $\mathrm{d}\Lambda
\left( t,a\right) =0\Rightarrow a=0$, where $0=ad_{t}$ for any $a\in 
\mathfrak{a}$, we can always identify $\mathfrak{a}$ with the linear span, 
\begin{equation*}
\sum \lambda _{i}\mathrm{d}\Lambda \left( t,a_{i}\right) =\mathrm{d}\Lambda
\left( t,\sum \lambda _{i}a_{i}\right) ,\quad \forall \lambda _{i}\in 
\mathbb{C},a_{i}\in \mathfrak{a},
\end{equation*}%
and consider it as a complex associative $\star $-algebra, having the death $%
d_{t}\in \mathfrak{a}$, a $\star $-invariant annihilator $\mathfrak{a}\cdot
d_{t}=\left\{ 0\right\} $ corresponding to $\mathrm{d}\Lambda \left( t,%
\mathfrak{a}\right) \mathrm{d}t=\left\{ 0\right\} $. The derivative $l$ of
the differential expectations $a\mapsto l\left( a\right) \mathrm{d}t$ with
respect to the Lebesgue measure $\mathrm{d}t$, called the It\^{o} algebra
state, is a linear positive $\star $-functional 
\begin{equation*}
l:\mathfrak{a}\rightarrow \mathbb{C},\quad l\left( a\cdot a^{\star }\right)
\geq 0,\quad l\left( a^{\star }\right) =\overline{l\left( a\right) },\quad
\forall a\in \mathfrak{a},
\end{equation*}%
normalized as $l\left( d_{t}\right) =1$ correspondingly to the determinism $%
\left\langle I\mathrm{d}t\right\rangle =\mathrm{d}t$ of $\mathrm{d}\Lambda
\left( t,d_{t}\right) $. We shall identify the It\^{o} algebra $\left( 
\mathrm{d}\Lambda \left( \mathfrak{a}\right) ,l\mathrm{d}t\right) $ and the
parametrizing algebra $\left( \mathfrak{a},l\right) $ and assume that it is
faithful in the sense that the $\star $-ideal 
\begin{equation}
\mathfrak{i}=\left\{ b\in \mathfrak{a}:l\left( b\right) =l\left( b\cdot
c\right) =l\left( a\cdot b\right) =l\left( a\cdot b\cdot c\right) =0\quad
\forall a,c\in \mathfrak{a}\right\}  \label{0.2}
\end{equation}%
is trivial, $\mathfrak{i}=\left\{ 0\right\} $, otherwise $\mathfrak{a}$
should be factorized with respect to this ideal. Note that the associativity
of the algebra $\mathfrak{a}$ as well as the possibility of its
noncommutativity is inherited from the associativity and noncommutativity of
the operator product $\Delta \Lambda \left( t,a\right) \Delta \Lambda \left(
t,b\right) $ on the pre-Hilbert space.

Now we can study the representations of the It\^{o} algebra $\left( 
\mathfrak{a},l\right) $. Because any It\^{o} algebra contains the Newton
nilpotent subalgebra $\left( \mathbb{C}d_{t},l\right) $, it has no identity
and cannnot be realized by operators in a Euclidean space even if it is
finite-dimensional $\star $-algebra. Thus we have to consider the operator
representations of $\mathfrak{a}$ in a pseudo-Euclidean space, and we shall
find such representations in the simplest one, in a complex Minkowski space.

Let $\mathbb{H}$ be a complex pseudo-Euclidean space with respect to a
separating indefinite metric $\left( \mathrm{x}|\mathrm{x}\right) $, and $%
\mathrm{h}\mathbf{\in }\mathbb{H}$ be a non-zero vector. We denote by $%
\mathcal{B}\left( \mathbb{H}\right) $ the $\dagger $-algebra of all
operators $\mathrm{A}:\mathbb{H}\rightarrow \mathbb{H}$ with $\mathrm{A}%
^{\dagger }\mathbb{H\subseteq H}$, where $\mathrm{A}^{\dagger }$ is defined
as the kernel of the Hermitian adjoint sesquilinear form $\left( \mathrm{x}|%
\mathrm{A}^{\dagger }\mathrm{x}\right) =\overline{\left( \mathrm{x}|\mathrm{%
Ax}\right) }$. A linear map $\mathrm{i}:\mathfrak{a}\rightarrow \mathcal{B}%
\left( \mathbb{H}\right) $ is a representation of the It\^{o} $\star $%
-algebra $\left( \mathfrak{a},l\right) $ on $\left( \mathbb{H},\mathrm{h}%
\right) $ if 
\begin{equation}
\mathrm{i}\left( a^{\star }\right) =\mathrm{i}\left( a\right) ^{\dagger
},\quad \mathrm{i}\left( a\cdot b\right) =\mathrm{i}\left( a\right) \mathrm{i%
}\left( b\right) ,\quad \left( \mathrm{h}\mathbf{|}\mathrm{i}\left( a\right) 
\mathrm{h}\right) =l\left( a\right) \quad \forall a,b\in \mathfrak{a}\text{.}
\label{0.5}
\end{equation}%
We can always assume that $\left( \mathrm{h}\mathbf{|}\mathrm{h}\right) =0$,
otherwise $\mathrm{h}$ should be replaced by the vector $\mathrm{h}_{+}=%
\mathrm{h}-\frac{1}{2}\left( \mathrm{h}|\mathrm{h}\right) \mathrm{h}_{-}$,
where $\mathrm{h}_{-}=\mathrm{i}\left( d_{t}\right) \mathrm{h}$, inducing
the same state 
\begin{eqnarray}
&&\left( \mathrm{h}_{+}|\mathrm{i}\left( a\right) \mathrm{h}_{+}\right)
\label{0.3} \\
&=&l\left( a\right) -\frac{1}{2}\left( \mathrm{h}|\mathrm{h}\right) \left( 
\mathrm{h}\mathbf{|}\mathrm{i}\left( d_{t}a+ad_{t}\emph{\,}-\emph{\,}\frac{1%
}{2}\left( \mathrm{h}|\mathrm{h}\right) d_{t}ad_{t}\right) \mathrm{h}\right)
\\
&=&l\left( a\right) \text{.}
\end{eqnarray}

\begin{proposition}
Every operator representation $\left( \mathbb{H},\mathrm{i},\mathrm{h}%
\right) $ of an It\^{o} algebra $\left( \mathfrak{a},l\right) $ is
equivalent to the triangular-matrix representation $\mathbf{i}=\left[ i_{\nu
}^{\mu }\right] _{\nu =-,\bullet ,+}^{\mu =-,\bullet ,+}$ with $i_{\nu
}^{\mu }\left( a\right) =0$ if $\mu =+$ or $\nu =-$ and $i_{+}^{-}\left(
a\right) =l\left( a\right) $ for all $a\in \mathfrak{a}$. Here $a_{\nu
}^{\mu }=i_{\nu }^{\mu }\left( a\right) $ are linear operators $\mathbb{H}%
_{\nu }\rightarrow \mathbb{H}_{\mu }$ on a pseudo-Hilbert (pre-Hilbert if
minimal) space $\mathbb{H}_{\bullet }$ and on $\mathbb{H}_{+}=\mathbb{C=H}%
_{-}$, having the adjoints $a_{\nu }^{\mu \dagger }:\mathbb{H}_{\mu
}\rightarrow \mathbb{H}_{\nu }$, which define the pseudo-Euclidean
involution $\mathbf{a}\mapsto \mathbf{a}^{\dagger }$ by $a_{-\nu }^{\star
\mu }=a_{-\mu }^{\nu \dagger }$, where $-\left( -,\bullet ,+\right) =\left(
+,\bullet ,-\right) $. Moreover, if the representation is minimal, then $%
\mathbb{H}_{\bullet }$ is a pre-Hilbert space and $i_{\nu }^{\mu }\left( 
\text{$d_{t}$}\right) =\delta _{-}^{\mu }\delta _{\nu }^{+}$.
\end{proposition}

\begin{proof}
In the matrix notation $i_{\nu }^{\mu }\left( a\right) =\mathrm{h}^{\mu }%
\mathrm{i}\left( a\right) \mathrm{h}_{\nu }$, where $\mathrm{h}^{-}=\mathrm{h%
}_{+}^{\dagger }$, $\mathrm{h}^{+}=\mathrm{h}_{-}^{\dagger }$ are defined by 
$\mathrm{h}^{\dagger }\mathrm{x}=\left( \mathrm{h}|\mathrm{x}\right) $ for
all $\mathrm{h}\mathbf{,}\mathrm{x}\in \mathbb{H}$, (\ref{0.3}) can be
written as $i_{+}^{-}\left( a\right) =l\left( a\right) $, and $%
i_{+}^{+}\left( a\right) =0=i_{-}^{-}\left( a\right) $ and $i_{-}^{+}\left(
a\right) =0$ as 
\begin{equation*}
\mathrm{i}\left( a\right) \mathrm{h}_{-}=\mathrm{i}\left( ad_{t}\right) 
\mathrm{h}=0=\mathrm{h}^{\dagger }\mathrm{i}\left( d_{t}a\right) =\mathrm{h}%
^{+}\mathrm{i}\left( a\right) \quad \forall a\in \mathfrak{a}\text{.}
\end{equation*}%
Moreover, due to the pseudo-orthogonality 
\begin{equation*}
\left( \mathrm{x}\mathbf{|}\mathrm{x}\right) =\xi _{-}\xi ^{-}+\left( 
\mathrm{x}_{\bullet }|\mathrm{x}_{\bullet }\right) +\xi _{+}\xi ^{+}\equiv
\left( \mathbf{x|x}\right) ,
\end{equation*}%
of the decomposition $\mathrm{x}=\xi ^{-}\mathrm{h}_{-}+\mathrm{x}_{\bullet
}+\xi ^{+}\mathrm{h}_{+}$, where $\xi ^{-}=\mathrm{h}^{-}\mathrm{x}=\bar{\xi}%
_{+}$, $\xi ^{+}=\mathrm{h}^{+}\mathrm{x}=\bar{\xi}_{-}$, $\mathbf{x}=\left(
\xi _{-},\mathrm{x}_{\bullet }^{\dagger },\xi _{+}\right) $, the
representation of the It\^{o} $\star $-algebra $\left( \mathfrak{a},l\right) 
$ is defined by the homomorphism $\mathbf{i}:a\mapsto \left[ i_{\nu }^{\mu
}\left( a\right) \right] $ into the space of triangular block-matrices $%
\mathbf{a}=\left[ a_{\nu }^{\mu }\right] _{\nu =-,\bullet ,+}^{\mu
=-,\bullet ,+}$ with $a_{\nu }^{\mu }=0$ if $\mu =+$ or $\nu =-$.

If the representation is minimal in the sense that $\mathbb{H}\mathbf{=}%
\mathrm{i}\left( \mathfrak{a}\right) \mathrm{h}$, and $\mathrm{h}$ has zero
length, it is pseudo-unitary equivalent to the triangular representation on
the complex Minkowski space $\mathbb{C}\oplus \mathbb{H}_{\bullet }\oplus 
\mathbb{C}$, as it can be easily seen in the basis $\mathrm{h}_{+}=\mathrm{h}
$, $\mathrm{h}_{-}=\mathrm{i}\left( d_{t}\right) \mathrm{h}$. Indeed, the
pseudo-orthogonal to the zero length vectors $\mathrm{h}_{-},\mathrm{h}_{+}$
space $\mathbb{H}_{\bullet }$ in this case is the pre-Hilbert space $\mathbb{%
H}_{\bullet }=\left\{ \mathrm{i}\left( a\right) \mathrm{h}:l\left( a\right)
=0\right\} $ as 
\begin{equation*}
\xi ^{+}=\mathrm{h}^{+}\mathrm{i}\left( a\right) \mathrm{h}=0,\quad \xi ^{-}=%
\mathrm{h}^{-}\mathrm{i}\left( a\right) \mathrm{h}=l\left( a\right) \quad
\forall a\in \mathfrak{a},
\end{equation*}%
and $\left( \mathrm{x}_{\bullet }|\mathrm{x}_{\bullet }\right) =l\left(
a^{\star }\cdot a\right) \geq 0$ for all $\mathrm{x}_{\bullet }=\mathrm{i}%
\left( a\right) \mathrm{h}-\xi ^{-}\mathrm{h}_{-}=\mathrm{i}\left( a-l\left(
a\right) d_{t}\right) \mathrm{h}$. Moreover, in the minimal representation $%
i_{\bullet }^{-}\left( d_{t}\right) =0=i_{+}^{\bullet }\left( d_{t}\right) $
and $i_{\bullet }^{\bullet }\left( d_{t}\right) =0$ as 
\begin{equation*}
\mathrm{i}\left( d_{t}\right) \mathrm{x}_{\bullet }=\mathrm{i}\left(
d_{t}\cdot a\right) \mathrm{h}=0=\mathrm{i}\left( a^{\star }\cdot
d_{t}\right) \mathrm{h}=\mathrm{x}_{\bullet }^{\dagger }\mathrm{i}\left(
d_{t}\right) \text{ \quad }\forall \mathrm{x}_{\bullet }\in \mathbb{H}%
_{\bullet }.
\end{equation*}%
Thus, the only nonzero matrix element of $\mathrm{i}\left( d_{t}\right) $ is 
$i_{+}^{-}\left( d_{t}\right) =1.$
\end{proof}

Note that the constructed equivalent matrix representation is also defined
as the right representation $\mathbf{x}\mapsto \mathbf{xa}$ on all
raw-vectors $\mathbf{x}=\left( \xi _{-},\xi ,\xi _{+}\right) $, $\xi \in 
\mathbb{H}_{\bullet }^{\dagger }$, into the dual space $\mathbb{H}^{*}=%
\mathbb{C\times H}_{\bullet }^{*}\mathbb{\times C\supseteq H}^{\dagger }$
with the invariant $\mathbb{H}^{\dagger }=\left\{ \mathrm{x}^{\dagger }:%
\mathrm{x}\in \mathbb{H}\right\} $ such that $a_{+}^{-}=\left( \mathrm{h}^{-}%
\mathbf{a|}\mathrm{h}^{-}\right) =l\left( a\right) $, where $\mathrm{h}%
^{-}=\left( 1,0,0\right) $. Moreover, in the Minkowski $\mathbb{C}\oplus 
\mathbb{H}_{\bullet }\oplus \mathbb{C}$ space it can be extended by a
continuity onto $\mathbb{C}\oplus \mathcal{H}\oplus \mathbb{C}$, where $%
\mathcal{H}$ is the closure of the pre-Hilbert space $\mathbb{H}_{\bullet }$
with respect to $\left\| \mathrm{x}_{\bullet }\right\| $ and all the
seminorms $\left\| i_{\bullet }^{\bullet }\left( a\right) \mathrm{x}%
_{\bullet }\right\| $, $a\in \mathfrak{a}$ simultaneously. We shall call
such representation closed if $\mathbb{H}$ is the minimal closed Minkowski
space $\mathbb{C}\oplus \mathcal{H}\oplus \mathbb{C}$, i.e. if $\mathcal{H}$
is the closure of the minimal pre-Hilbert space $\mathbb{H}_{\bullet
}=i_{+}^{\bullet }\left( \mathfrak{a}\right) $.

\begin{theorem}
Every It\^{o} $\star $-algebra $\left( \mathfrak{a},l\right) $ can be
canonically realized as the triangular matrix subalgebra in a complex
Minkowski space. Moreover, every minimal closed representation is equivalent
to the canonical one.
\end{theorem}

\begin{proof}
Now we shall construct a faithful canonical operator representation for any
It\^{o} algebra $\left( \mathfrak{a},l\right) $. The functional $l$ defines
for each $a\in \mathfrak{a}$ the canonical quadruple 
\begin{equation}
a_{\bullet }^{\bullet }=i\left( a\right) ,\quad a_{+}^{\bullet }=k\left(
a\right) ,\quad a_{\bullet }^{-}=k^{\dagger }\left( a\right) ,\quad
a_{+}^{-}=l\left( a\right) ,  \label{0.6}
\end{equation}
where $i\left( a\right) =i\left( a^{\star }\right) ^{\dagger }$ is the GNS
representation $k\left( a\cdot b\right) =i\left( a\right) k\left( b\right) $
of $\mathfrak{a}$ on the pre-Hilbert space $\mathbb{H}_{\bullet }=\left\{
k\left( b\right) :b\in \mathfrak{a}\right\} $ of the Kolmogorov
decomposition $l\left( a\cdot b\right) =k^{\dagger }\left( a\right) k\left(
b\right) $, and $k^{\dagger }\left( a\right) =k\left( a^{\star }\right)
^{\dagger }$. Because the operators $i\left( a\right) $ are continuous on $%
\mathbb{H}_{\bullet }$ w.r.t. the topology induced by $\left\| k\left(
b\right) \right\| $ and all the seminorms $\left\| k\left( a\cdot b\right)
\right\| $, $a\in \mathfrak{a}$ on $\mathbb{H}_{\bullet }$, can define the
representation $i$ on the clousure $\mathcal{H}=\overline{k\left( \mathfrak{a%
}\right) }$ of $\mathbb{H}_{\bullet }$ w.r.t. these seminorms. The obtained
quadrupole representation $\boldsymbol{i}:a\mapsto \boldsymbol{a}=\left(
a_\nu ^\mu \right) _{\nu =+,\bullet }^{\mu =-,\bullet }$ of $\mathfrak{a}$
is multiplicative, $\boldsymbol{i}\left( a\cdot b\right) =\left( a_{\bullet
}^\mu b_\nu ^{\bullet }\right) _{\nu =+,\bullet }^{\mu =-,\bullet }$ with
respect to the product given by the convolution of the components $a_\nu $
and $b^\mu $ over the common index values $\mu =\bullet =\nu $: 
\begin{eqnarray*}
i\left( a\right) i\left( b\right) &=&i\left( a\cdot b\right) ,\quad
k^{\dagger }\left( a\right) i\left( b\right) =k^{\dagger }\left( a\cdot
b\right) \\
i\left( a\right) k\left( b\right) &=&k\left( a\cdot b\right) ,\quad
k^{\dagger }\left( a\right) k\left( b\right) =l\left( a\cdot b\right) \text{.%
}
\end{eqnarray*}
It is faithful because of the triviality of the ideal (\ref{0.2}). Now we
can use the convenience $a_{-}^\mu =0=a_\nu ^{+}$ of the tensor notations (%
\ref{0.6}), extending the quadruples $\boldsymbol{a}=\boldsymbol{i}\left(
a\right) $ to the triangular matrices $\mathbf{a}=\left[ a_\nu ^\mu \right]
_{\nu =-,\bullet ,+}^{\mu =-,\bullet ,+}$, in which (\ref{0.4}) is simply
given by $\mathbf{i}\left( a\cdot b\right) =\mathbf{ab}$ in terms of the
usual product of the matrices $\mathbf{a}=\mathbf{i}\left( a\right) $ and $%
\mathbf{b}=\mathbf{i}\left( b\right) $. However the involution $a\mapsto
a^{\star }$, which is given by the Hermitian conjugation $\boldsymbol{i}%
\left( a^{\star }\right) =\left( a_{-\mu }^{-\nu \dagger }\right) _{\nu
=+,\bullet }^{\mu =-,\bullet }$ of the quadruples $\boldsymbol{a}$, where $%
-(-)=+$, $-\bullet =\bullet $, $-(+)=-$, is represented by the adjoint
matrix $\mathbf{a}^{\dagger }=\mathrm{G}\mathbf{a}^{*}\mathrm{G}$ w.r.t. the
pseudo-Euclidean (complex Minkowski) metric tensor $\mathrm{G}=\left[ \delta
_{-\nu }^\mu \right] _{\nu =-,\bullet ,+}^{\mu =-,\bullet ,+}$. Thus, we
have constructed the canonical representation 
\begin{equation*}
\mathbf{i}\left( a\right) =\left[ 
\begin{array}{lll}
0 & k^{\dagger }\left( a\right) & l\left( a\right) \\ 
0 & i\left( a\right) & k\left( a\right) \\ 
0 & 0 & 0%
\end{array}
\right] ,\quad \mathbf{i}\left( a^{\star }\right) =\mathrm{G}\mathbf{i}%
\left( a\right) ^{*}\mathrm{G}\mathbf{,\quad }\mathrm{G}=\left[ 
\begin{array}{lll}
0 & 0 & 1 \\ 
0 & I & 0 \\ 
1 & 0 & 0%
\end{array}
\right]
\end{equation*}
in the Minkowski space $\mathbb{C\oplus }\mathcal{H}\oplus \mathbb{C}$ with $%
\mathcal{H}=\overline{k\left( \mathfrak{a}\right) }$ and $\mathrm{h}%
^{-}=\left( 1,0,0\right) $.

The second part of the Theorem follows from the fact that every minimal
representation space is the Minkowski one. All Minkowski spaces of the
minimal closed representations are pseudo-unitary equivalent because all
minimal closed pre-Hilbert spaces containing $i_{+}^{\bullet }\left( 
\mathfrak{a}\right) $ are unitary equivalent.
\end{proof}

\begin{definition}
Let $\mathcal{H}$ be a pre-Hilbert space, and $\mathfrak{b}\left( \mathcal{H}%
\right) $ be the associated $\star $-algebra of all quadrupoles $A=\left(
a_\nu ^\mu \right) _{\nu =+,\bullet }^{\mu =-,\bullet }$, where $a_\nu ^\mu $
are linear operators $\mathbb{H}_\nu \rightarrow \mathbb{H}_\mu $ with $%
\mathbb{H}_{\bullet }=\mathcal{H}$, $\mathbb{H}_{+}=\mathbb{C=H}_{-}$,
having the adjoints $a_\nu ^{\mu \dagger }:\mathbb{H}_\mu \rightarrow 
\mathbb{H}_\nu $, with the product and the involution 
\begin{equation}
A\cdot B=\left( a_{\bullet }^\mu b_\nu ^{\bullet }\right) _{\nu =+,\bullet
}^{\mu =-,\bullet },\quad A^{\star }=\left( a_{-\mu }^{-\nu \dagger }\right)
_{\nu =+,\bullet }^{\mu =-,\bullet }\text{.}  \label{0.4}
\end{equation}
It is an It\^{o} algebra with respect to $l\left( \boldsymbol{a}\right)
=a_{+}^{-}$ and the death $D_{t}=\left( \delta _{-}^\mu \delta _\nu
^{+}\right) _{\nu =+,\bullet }^{\mu =-,\bullet }=D_{t}^{\star }$, $
A\cdot D_{t}=0$, $\forall A\in \mathfrak{b}\left( \mathcal{H}%
\right) $, called the HP (Hudson-Parthasarathy) algebra associated with the
space $\mathcal{H}$. The fundamental representation of an It\^{o} algebra $%
\left( \mathfrak{a},l\right) $ is given by the constructed homomorphism $%
\boldsymbol{i}:\mathfrak{a}\rightarrow \mathfrak{b}\left( \mathcal{H}\right) 
$ 
\begin{equation}
\boldsymbol{i}\left( a\right) =\left( 
\begin{tabular}{ll}
$l\left( a\right) $ & $k^{\dagger }\left( a\right) $ \\ 
$k\left( a\right) $ & $i\left( a\right) $%
\end{tabular}
\right) ,\quad \,\boldsymbol{i}\left( a^{\star }\right) =\left( 
\begin{tabular}{ll}
$l\left( a^{\star }\right) $ & $k\left( a\right) ^{\dagger }$ \\ 
$k\left( a^{\star }\right) $ & $i\left( a\right) ^{\dagger }$%
\end{tabular}
\right)  \label{0.7}
\end{equation}
$\,\boldsymbol{i}\left( a\cdot b\right) =\boldsymbol{i}\left( a\right) \cdot 
\boldsymbol{i}\left( b\right) $, into the HP algebra, associated with the
space $\mathcal{H}$ of its canonical representation.
\end{definition}

Note that because the It\^{o} algebra is assumed faithful in the sense of
the triviality of the ideal (\ref{0.2}), the fundamental representation, and
so the canonical representation (\ref{0.7}), is also faithful. It proves in
this case that $\Lambda \left( t,a\right) =a_\nu ^\mu \Lambda _\mu ^\nu
\left( t\right) $, where 
\begin{equation}
a_\nu ^\mu \Lambda _\mu ^\nu \left( t\right) =a_{\bullet }^{\bullet }\Lambda
_{\bullet }^{\bullet }\left( t\right) +a_{+}^{\bullet }\Lambda _{\bullet
}^{+}\left( t\right) +a_{\bullet }^{-}\Lambda _{-}^{\bullet }\left( t\right)
+a_{+}^{-}\Lambda _{-}^{+}\left( t\right) ,  \label{0.8}
\end{equation}
is the canonical decomposition of $\Lambda $ into the exchange $\Lambda
_{\bullet }^{\bullet }$, creation $\Lambda _{\bullet }^{+}$, annihilation $%
\Lambda _{-}^{\bullet }$ and preservation (time) $\Lambda _{-}^{+}=t\mathrm{I%
}$ operator-valued processes of the HP quantum stochastic calculus, having
the mean values $\left\langle \Lambda _\mu ^\nu \left( t\right)
\right\rangle =t\delta _{+}^\nu \delta _\mu ^{-}$. This was already noted in 
\cite{3,4} that any (classical or quantum) stochastic noise described by a
process $t\in \mathbb{R}_{+}\mapsto \Lambda \left( t,a\right) ,a\in 
\mathfrak{a}$ with independent increments $\mathrm{d}\Lambda \left(
t,a\right) =\Lambda \left( t+\mathrm{d}t,a\right) -\Lambda \left( t,a\right) 
$ forming an It\^{o} $\dagger $-algebra, can be represented in the Fock
space $\mathfrak{F}$ over the space of $\mathcal{H}$-valued
square-integrable functions on $\mathbb{R}_{+}$, with the vacuum vector
state.

\section{Two basic It\^{o} B*-algebras}

Here we consider two extreme cases of Banach It\^{o} algebras as closed
sub-algebras of the vaccum HP-algebra $\mathfrak{b}\left( \mathcal{H}\right) 
$ associated with a Hilbert space $\mathcal{H}$. The first case correspondes
to a pure state $l$ on $\mathfrak{a}$ as it is in the case of a quantum
noise of zero temperature, and the second case corresponds to a completely
mixed $l$ as in the case of a quantum noise of a finite temperature.

\subsection{Vacuum noise B*-algebra}

Let $\mathcal{H}$ be a Hilbert space of ket-vectors $\zeta $ with scalar
product $\left( \zeta |\zeta \right) \equiv \zeta ^{\dagger }\zeta $ and $%
\mathcal{A}\subseteq \mathcal{B}\left( \mathcal{H}\right) $ be a C*-algebra,
represented on $\mathcal{H}$ by the operators $\mathcal{A}\ni A:\zeta
\mapsto A\zeta $ with $\left( A^{\dagger }\zeta |\xi \right) =\left( \zeta
|A\xi \right) $, $\xi \in \mathcal{H}$. We denote by $\mathcal{H}^{\dagger }$
the dual Hilbert space of bra-vectors $\eta =\zeta ^{\dagger }$, $\zeta \in 
\mathcal{H}$ with the scalar product $\left( \eta |\xi ^{\dagger }\right)
=\eta \xi =\left( \eta ^{\dagger }|\xi \right) $ given by inverting
anti-linear isomorphism $\mathcal{H}^{\dagger }\ni \eta \mapsto \eta
^{\dagger }\in \mathcal{H}$, and the dual representation of $\mathcal{A}$ as
the right representation $A^{\prime }:\eta \mapsto \eta A$, $\eta \in 
\mathcal{H}^{\dagger }$, given by $\left( \eta A\right) \zeta =\eta \left(
A\zeta \right) $ such that $\left( \eta A^{\dagger }|\eta \right) =\left(
\eta |\eta A\right) $ on $\mathcal{H}^{\dagger }$. Then the direct sum $%
\mathcal{K}=\mathcal{H}\oplus \mathcal{H}^{\dagger }$ of $\xi =\zeta \oplus
\eta $ becomes a two-sided $\mathcal{A}$-module 
\begin{equation}
A\left( \zeta \oplus \eta \right) =A\zeta ,\quad \left( \zeta \oplus \eta
\right) A=\eta A,\quad \forall \zeta \in \mathcal{H},\eta \in \mathcal{H}%
^{\dagger },  \label{1.1}
\end{equation}%
with the flip-involution $\xi ^{\star }=\eta ^{\dagger }\oplus \zeta
^{\dagger }$ and two scalar products 
\begin{equation}
\left\langle \zeta \oplus \eta ^{\prime }|\zeta ^{\prime }\oplus \eta
\right\rangle _{+}=\left( \zeta |\zeta ^{\prime }\right) ,\quad \left\langle
\zeta \oplus \eta ^{\prime }|\zeta ^{\prime }\oplus \eta \right\rangle
^{-}=\left( \eta ^{\prime }|\eta \right) .  \label{1.2}
\end{equation}

The space $\mathfrak{a}=\mathbb{C}\oplus \mathcal{K}\oplus \mathcal{A}$ of
triples $a=\left( \alpha ,\xi ,A\right) $ becomes an It\^{o} $\star $%
-algebra with respect to the non-commutative product 
\begin{equation}
a^{\star }\cdot a=\left( \left\langle \xi |\xi \right\rangle _{+},\xi
^{\star }A+A^{\dagger }\xi ,A^{\dagger }A\right) ,\quad a\cdot a^{\star
}=\left( \left\langle \xi |\xi \right\rangle ^{-},A\xi ^{\star }+\xi
A^{\dagger },AA^{\dagger }\right) ,  \label{1.3}
\end{equation}%
where $\left( \alpha ,\xi ,A\right) ^{\star }=\left( \bar{\alpha},\xi
^{\star },A^{\dagger }\right) $, with death $d_{t}=\left( 1,0,0\right) $ and 
$l\left( \alpha ,\xi ,A\right) =\alpha $. Obviously $a^{\star }\cdot a\neq
a\cdot a^{\star }$ if $\left\Vert \xi \right\Vert _{+}=\left\Vert \zeta
\right\Vert \neq $ $\left\Vert \eta \right\Vert =\left\Vert \xi \right\Vert
^{-}$ even if the operator algebra $\mathcal{A}$ is commutative, $A^{\dagger
}A=AA^{\dagger }$. It is separated by four semi-norms 
\begin{equation}
\left\Vert a\right\Vert =\left\Vert A\right\Vert ,\,\left\Vert a\right\Vert
_{+}=\left\Vert \zeta \right\Vert ,\,\left\Vert a\right\Vert ^{-}=\left\Vert
\eta \right\Vert ,\,\left\Vert a\right\Vert _{+}^{-}=\left\vert \alpha
\right\vert ,  \label{1.4}
\end{equation}%
and is jointly complete as $a=\left( \alpha ,\zeta \oplus \eta ,A\right) \in 
\mathfrak{a}$ have independent components from the Banach spaces $\mathbb{C}$%
, $\mathcal{H}$, $\mathcal{H}^{\dagger }$and $\mathcal{A}$.

We shall call such Banach It\^{o} algebra the vacuum algebra as $l\left(
a^{\star }\cdot a\right) =0$ for any $a\in \mathfrak{a}$ with $\xi \in 
\mathcal{H}^{\dagger }$ (it is Hudson-Parthasarathy algebra $\mathfrak{a}=%
\mathfrak{b}\left( \mathcal{H}\right) $ if $\mathcal{A}=\mathcal{B}\left( 
\mathcal{H}\right) $). Every closed It\^{o} subalgebra $\mathfrak{a}%
\subseteq \mathfrak{b}\left( \mathcal{H}\right) $ of the HP algebra $%
\mathfrak{b}\left( \mathcal{H}\right) $ equipped with four norms (\ref{1.4})
on a Hilbert space $\mathcal{H}$ is called the operator It\^{o} B*-algebra.

If the algebra $\mathcal{A}$ is completely degenerated on $\mathcal{H}$, $%
\mathcal{A}=\left\{ 0\right\} $, the It\^{o} algebra $\mathfrak{a}$ is
nilpotent of second order, and contains only the two-dimensional subalgebras
of Wiener type $\mathfrak{b}=\mathbb{C}\oplus \mathbb{C}\oplus \left\{
0\right\} $ generated by an $a=\left( \alpha ,\zeta \oplus \eta ,0\right) $
with $\left\| \zeta \right\| =\left\| \eta \right\| $. Every closed It\^{o}
subalgebra $\mathfrak{b}\subseteq \mathfrak{a}$ of the HP B*-algebra $%
\mathfrak{a}=\mathfrak{b}\left( \mathcal{H}\right) $ is called the B*-It\^{o}
algebra of a vacuum Brownian motion if it is defined by a $\star $-invariant
direct sum $\mathcal{G}=\mathcal{G}_{+}\oplus \mathcal{G}^{-}\subseteq 
\mathcal{K}$ given by a Hilbert subspace $\mathcal{G}_{+}\subseteq \mathcal{H%
}$, $\mathcal{G}^{-}=\mathcal{G}_{+}^{\dagger }$ and $\mathcal{A}=\left\{
0\right\} $.

In the case $I\in \mathcal{A}$ the algebra $\mathcal{A}$ is not degenerated
and contains also the vacuum Poisson subalgebra $\mathbb{C}\oplus \left\{
0\right\} \oplus \mathbb{C}I$ of the total quantum number on $\mathcal{H}$,
and other Poisson two-dimensional subalgebras, generated by $a=\left( \alpha
,\zeta \oplus \eta ,I\right) $ with $\eta =e^{i\theta }\zeta ^{\dagger }$.
We shall call a closed It\^{o} subalgebra $\mathfrak{c}\subseteq \mathfrak{a}
$ of the HP B*-algebra $\mathfrak{a}=\mathfrak{b}\left( \mathcal{H}\right) $
the B*-algebra of a vacuum L\'{e}vy motion if it is given by a direct sum $%
\mathcal{E}=\mathcal{E}_{+}\oplus \mathcal{E}^{-}\subseteq \mathcal{K}$ with 
$\mathcal{E}^{\_}=\mathcal{E}_{+}^{\dagger }$ and a $\dagger $-subalgebra $%
\mathcal{A}\subseteq \mathcal{B}\left( \mathcal{H}\right) $ nondegenerated
on the subspace $\mathcal{E}_{+}\subseteq \mathcal{H}$.

\begin{theorem}
Every vacuum B*-algebra can be decomposed into an orthogonal sum $\mathfrak{a%
}=\mathfrak{b}+\mathfrak{c}$, $\mathfrak{b}\cdot \mathfrak{c}=\left\{
0\right\} $ of the Brownian vacuum B*-algebra $\mathfrak{b}$ and the L\'{e}%
vy vacuum B*-algebra $\mathfrak{c}$. This decomposition is unique on the
zero mean kernel $\mathfrak{x}=\left\{ x\in \mathfrak{a}:l\left( x\right)
=0\right\} $.
\end{theorem}

\begin{proof}
This decomposition is uniquely defined for all $a=\left( \alpha ,\xi
,A\right) $ by $a=\alpha d_{t}+y+z$, with $y=\left( 0,\eta ,0\right) $, $%
z=\left( 0,\zeta ,A\right) $, $\eta =P\xi \oplus \xi P\in \mathcal{G}$, $%
\zeta =\xi -\eta \in \mathcal{E}$, where $P=P^{\dagger }$ is the maximal
projector in $\mathcal{H}$, for which $\mathcal{A}P=\left\{ 0\right\} $, $%
\mathcal{G}_{+}=P\mathcal{H}$, and $\mathcal{E}_{+}\mathcal{=G}_{+}^{\perp }$%
.
\end{proof}

\subsection{Thermal noise B*-algebra}

Let $\mathcal{D}$ be a left Tomita $\star $-algebra \cite{5} with respect to
a Hilbert norm $\left\Vert \xi \right\Vert _{+}=0\Rightarrow \xi =0$, and
thus a right pre-Hilbert $\star $-algebra with respect to $\left\Vert \xi
\right\Vert ^{-}=\left\Vert \xi ^{\star }\right\Vert _{+}$. This means that $%
\mathcal{D}$ is a complex pre-Hilbert space with continuous left (right)
multiplications $C:\zeta \mapsto \xi \zeta $ ($C^{\prime }:\eta \mapsto \eta
\xi $) w.r.t. $\left\Vert \cdot \right\Vert _{+}$ (w.r.t. $\left\Vert \cdot
\right\Vert ^{-}$) of the elements $\zeta ,\eta \in \mathcal{D}$
respectively, defined by an associative product in $\mathcal{D}$, and the
involution $\mathcal{D}\ni \xi \mapsto \xi ^{\star }\in \mathcal{D}$ such
that 
\begin{equation}
\left\langle \eta \zeta ^{\star }|\xi \right\rangle ^{-}=\left\langle \eta
|\xi \zeta \right\rangle ^{-},\quad \left\langle \eta ^{\star }\zeta |\xi
\right\rangle _{+}=\left\langle \zeta |\eta \xi \right\rangle _{+}\quad
\forall \xi ,\zeta ,\eta \in \mathcal{D},  \label{2.1}
\end{equation}%
\begin{equation}
\left\langle \eta |\xi ^{\star }\right\rangle ^{-}=\left\langle \xi |\eta
^{\sharp }\right\rangle ^{-},\quad \left\langle \zeta |\xi ^{\star
}\right\rangle _{+}=\left\langle \xi |\zeta ^{\flat }\right\rangle _{+}\quad
\forall \eta \in \mathcal{D}^{-},\zeta \in \mathcal{D}_{+}.  \label{2.2}
\end{equation}%
Here $\left\langle \eta |\xi ^{\star }\right\rangle ^{-}=\left\langle \eta
^{\star }|\xi \right\rangle _{+}$ is the right scalar product, $\mathcal{D}%
_{+}=\mathcal{D}_{+}^{\flat }$ is a dense domain for the left adjoint
involution $\zeta \mapsto \zeta ^{\flat }$, $\zeta ^{\flat \flat }=\zeta $,
and $\mathcal{D}^{-}=\mathcal{D}_{+}^{\star }$ is the invariant domain for
the right adjoint involution $\eta \mapsto \eta ^{\sharp }$, $\left( \eta
^{\sharp }\eta \right) ^{\sharp }=\eta ^{\sharp }\eta $ such that $\zeta
^{\flat \star }=\zeta ^{\star \sharp }$, $\eta ^{\sharp \star }=\eta ^{\star
\flat }$.

Since the adjoint operators $C^{\dagger }\zeta =\xi ^{\star }\zeta $, $\eta
C^{\dagger }=\eta \xi ^{\star }$ are also given by the multiplications, they
are bounded: 
\begin{equation}
\left\| \xi \right\| =\sup \left\{ \left\| \xi \zeta \right\| _{+}:\left\|
\zeta \right\| _{+}\leq 1\right\} =\sup \left\{ \left\| \eta \xi \right\|
^{-}:\left\| \eta \right\| ^{-}\leq 1\right\} <\infty .  \label{2.3}
\end{equation}
Note that we do not require the sub-space $\mathcal{DD}\subseteq \mathcal{D}$
of all products $\left\{ \eta \zeta :\eta ,\zeta \in \mathcal{D}\right\} $
to be dense in $\mathcal{D}$ w.r.t. any of two Hilbert norms on $\mathcal{D}$%
, but it is always dense w.r.t. the operator semi-norm (\ref{2.3}) on $%
\mathcal{D}$. Hence the operator $\dagger $-algebra $\mathcal{C}=\left\{ C:%
\mathcal{D}\ni \zeta \mapsto \xi \zeta |\xi \in \mathcal{D}\right\} $ w.r.t.
the left scalar product, which is also represented on the $\mathcal{D}\ni
\eta $ equipped with $\left\langle \cdot |\cdot \right\rangle ^{-}$ by the
right multiplications $\eta C=\eta \xi $, $\xi \in \mathcal{D}$, can be
degenerated on $\mathcal{D}$.

Thus the direct sum $\mathfrak{a}=\mathbb{C}\oplus \mathcal{D}$ of pairs $%
a=\left( \alpha ,\xi \right) $ becomes an It\^{o} $\star $-algebra with the
product 
\begin{equation}
a^{\star }\cdot a=\left( \left\langle \xi |\xi \right\rangle _{+},\xi
^{\star }\xi \right) ,\quad a\cdot a^{\star }=\left( \left\langle \xi |\xi
\right\rangle ^{-},\xi \xi ^{\star }\right) ,  \label{2.4}
\end{equation}%
where $\left( \alpha ,\xi \right) ^{\star }=\left( \bar{\alpha},\xi ^{\star
}\right) $, with death $d_{t}=\left( 1,0\right) $ and $l\left( \alpha ,\xi
\right) =\alpha $. Obviously $a^{\star }a\neq aa^{\star }$ if the involution 
$a\mapsto a^{\star }$ is not isometric w.r.t. any of two Hilbert norms even
if the algebra $\mathcal{D}$ is commutative. It is a Banach algebra
separated by the semi-norms 
\begin{equation}
\left\Vert a\right\Vert =\left\Vert \xi \right\Vert ,\left\Vert a\right\Vert
^{-}=\left\Vert \xi \right\Vert ^{-},\left\Vert a\right\Vert _{+}=\left\Vert
\xi \right\Vert _{+},\left\Vert a\right\Vert _{+}^{-}=\left\vert \alpha
\right\vert .  \label{2.5}
\end{equation}%
if its normed $\star $-algebra $\mathcal{D}$ is complete jointly w.r.t. to
the first three norms.

We shall call such complete It\^{o} algebra the thermal B*-algebra as $%
l\left( a^{\star }a\right) =\left\| \xi \right\| _{+}^2\neq 0$ for any $a\in 
\mathfrak{a}$ with $\xi \neq 0$. If $\zeta \eta =0$ for all $\zeta ,\eta \in 
\mathcal{D}$, it is the It\^{o} B*-algebra of thermal Brownian motion. A
thermal B*-subalgebra $\mathfrak{b}\subseteq \mathfrak{a}$ with such trivial
product is given by any involutive pre-Hilbert $\star $-invariant two-normed
subspace $\mathcal{G}\subseteq \mathcal{D}$ which is closed w.r.t. the
Hilbert sum $\left\langle \eta |\zeta \right\rangle ^{-}+\left\langle \zeta
|\eta \right\rangle _{+}$. We shall call such Brownian algebra $\mathfrak{b}=%
\mathbb{C\oplus }\mathcal{G}$ the quantum (if $\left\| \cdot \right\|
_{+}\neq $ $\left\| \cdot \right\| ^{-}$) Wiener B*-algebra associated with
the space $\mathcal{G}$.

In the opposite case, if $\mathcal{DD}=\left\{ \zeta \eta :\zeta ,\eta \in 
\mathcal{D}\right\} $ is dense in $\mathcal{D}$, it has nondegenerated
operator representation $\mathcal{C}$ on $\mathcal{D}$. Any closed
involutive sub-algebra $\mathcal{E}\subseteq \mathcal{D}$ which is
non-degenerated on $\mathcal{E}$ defines an It\^{o} B*-algebra $\mathfrak{c}=%
\mathbb{C}\oplus \mathcal{E}$ of thermal L\'{e}vy motion. We shall call such
It\^{o} algebra the quantum (if $\mathcal{E}$ is non-commutative) Poisson
B*-algebra.

\begin{theorem}
Every thermal It\^{o} B*-algebra is an orthogonal sum $\mathfrak{a}=%
\mathfrak{b}+\mathfrak{c}$, $\mathfrak{bc}=\left\{ 0\right\} $ of the Wiener
B*-algebra $\mathfrak{b}=\mathfrak{b}^{\star }$ and the Poisson B*-algebra $%
\mathfrak{c}=\mathfrak{c}^{\star }$. This decomposition is unique on the
zero mean kernel $\mathfrak{x}=l^{-1}\left( 0\right) $.
\end{theorem}

\begin{proof}
The orthogonal decomposition $a=\alpha d_{t}+b+c$ for all $a=\left( \alpha
,\xi \right) \in \mathfrak{a}$, uniquely given by the decomposition $\xi
=\eta +\zeta $, where $\eta $ is the orthogonal projection of $\xi $ onto
the orthogonal complement $\mathcal{G}$ of $\mathcal{DD}$ w.r.t. any of two
scalar products in $\mathcal{D}$, and $\zeta =\xi -\eta $ .

Indeed, if $\xi \in \mathcal{D}$ is left orthogonal to $\mathcal{DD}$, then
it is also right orthogonal to $\mathcal{DD}$ and vice versa: 
\begin{eqnarray*}
\left\langle \eta \zeta ^{\star }|\xi \right\rangle ^{-} &=&\left\langle
\zeta \eta ^{\star }|\xi ^{\star }\right\rangle _{+}=\left\langle \xi |\eta
^{\sharp \star }\zeta ^{\flat }\right\rangle _{+}=0,\quad \forall \eta \in 
\mathcal{D}^{-},\zeta \in \mathcal{D}_{+}, \\
\left\langle \eta ^{\star }\zeta |\xi \right\rangle _{+} &=&\left\langle
\zeta ^{\star }\eta |\xi ^{\star }\right\rangle ^{-}=\left\langle \xi |\eta
^{\sharp }\zeta ^{\flat \star }\right\rangle ^{-}=0,\quad \forall \eta \in 
\mathcal{D}^{-},\zeta \in \mathcal{D}_{+},
\end{eqnarray*}
and $\xi ^{\sharp },\xi ^{\flat }$ are also orthogonal to $\mathcal{DD}$: 
\begin{equation*}
\left\langle \eta \zeta ^{\star }|\xi ^{\sharp }\right\rangle
^{-}=\left\langle \xi |\zeta \eta ^{\star }\right\rangle ^{-}=0,\quad
\left\langle \eta ^{\star }\zeta |\xi ^{\flat }\right\rangle
_{+}=\left\langle \xi |\zeta ^{\star }\eta \right\rangle _{+}=0
\end{equation*}
From these and (\ref{2.1}) equations it follows that $\eta \xi =0=\xi \zeta $
for all $\zeta ,\eta \in \mathcal{D}$ if $\xi $ is (right or left)
orthogonal to $\mathcal{DD}$, and so $\left\| \xi \right\| =0$ for such $\xi 
$ and vice versa. Thus the maximal orthogonal subspace to $\mathcal{DD}$ is
the $\star $-invariant space $\mathcal{G}=\left\{ \xi \in \mathcal{D}%
:\left\| \xi \right\| =0\right\} $, which is complete w.r.t. the two norms $%
\left\| \cdot \right\| ^{-},\left\| \cdot \right\| _{+}$, i.e. is a closed
subspace of the pre-Hilbert space $\mathcal{D}$ with the $\star $-invariant
scalar product 
\begin{equation*}
\left\langle \zeta ^{\star }|\eta \right\rangle =\left\langle \zeta ^{\star
}|\eta \right\rangle _{+}+\left\langle \eta |\zeta ^{\star }\right\rangle
^{-}=\left\langle \zeta |\eta ^{\star }\right\rangle ^{-}+\left\langle \eta
^{\star }|\zeta \right\rangle _{+}=\left\langle \eta ^{\star }|\zeta
\right\rangle .
\end{equation*}
Denoting by $P$ the right orthogonal projector in $\mathcal{D}$ onto the
orthogonal complement to $\mathcal{DD}$, we obtain $P\xi =\eta \in \mathcal{G%
}$, $P\xi ^{\star }=\eta ^{\star }$ as 
\begin{equation*}
\left\langle \eta ^{\star }|P\xi \right\rangle _{+}=\left\langle P\xi |\eta
^{\sharp }\right\rangle ^{-}=\left\langle \xi |\eta ^{\sharp }\right\rangle
^{-}=\left\langle \eta |\xi ^{\star }\right\rangle ^{-}=\left\langle \eta
|P\xi ^{\star }\right\rangle ^{-},\,\eta \in \mathcal{G}^{-},\xi \in 
\mathcal{D}
\end{equation*}
and thus $P\xi $ is also the left projection, and so the symmetric
projection $\left\langle \eta |P\xi \right\rangle =\left\langle \eta |\xi
\right\rangle $ of $\xi \in \mathcal{D}$ onto $\mathcal{G}\ni \eta $. So $%
\zeta =\xi -\eta \in \mathcal{D}$ is in the left (right) closure $\mathcal{E}%
\subseteq \mathcal{D}$ of $\mathcal{DD}$ .
\end{proof}

\section{Decomposition of It\^{o} B*-algebras}

Now we shall study the general uniformly bounded infinite dimensional Ito $%
\star $-algeblas, unifying the considered two basic cases.

Let $\mathfrak{a}$ be an associative infinite-dimensional complex algebra
with involution $b^{\star }=a\in \mathfrak{a},\forall b=a^{\star }$ which is
defined by the properties

\begin{equation*}
\left( b\cdot b^{\star }\right) ^{\star }=b\cdot b^{\star },\quad \left(
\sum \lambda _ib_i\right) ^{\star }=\sum \bar{\lambda}_ib_i^{\star },\quad
\forall b_i\in \mathfrak{a},\lambda _i\in \mathbb{C}.
\end{equation*}
We shall suppose that this algebra is a normed space with respect to four
semi-norms $\left\| \cdot \right\| _\nu ^\mu $, indexed as 
\begin{equation}
\left\| \cdot \right\| _{\bullet }^{\bullet }\equiv \left\| \cdot \right\|
,\quad \left\| \cdot \right\| _{+}^{\bullet }\equiv \left\| \cdot \right\|
_{+},\quad \left\| \cdot \right\| _{\bullet }^{-}\equiv \left\| \cdot
\right\| ^{-},\quad \left\| \cdot \right\| _{+}^{-}  \label{3.1}
\end{equation}
by $\mu =-,\bullet $, $\nu =+,\bullet $, satisfying the following conditions

\begin{equation}
\left\| b^{\star }\right\| =\left\| b\right\| ,\quad \left\| b^{\star
}\right\| _{+}=\left\| b\right\| ^{-},\quad \left\| b^{\star }\right\|
_{+}^{-}=\left\| b\right\| _{+}^{-},  \label{3.2}
\end{equation}
\begin{equation}
\left( \left\| a\cdot c\right\| _\nu ^\mu \leq \left\| a\right\| _{\bullet
}^\mu \left\| c\right\| _\nu ^{\bullet },\right) _{\nu =+,\bullet }^{\mu
=-,\bullet }\quad \forall \,a,b,c\in \mathfrak{a}.  \label{3.3}
\end{equation}
Thus the semi-norms (\ref{3.1}) separate $\mathfrak{a}$ in the sense 
\begin{equation*}
\left\| a\right\| =\left\| a\right\| _{+}=\left\| a\right\| ^{-}=\left\|
a\right\| _{+}^{-}=0\quad \Rightarrow \quad a=0,
\end{equation*}
and the product $\left( a,c\right) \mapsto ac$ with involution $\star $ is
uniformly continuous in the induced topology due to (\ref{3.3}).

If $\mathfrak{a}$ is a $\star $-algebra equipped with a linear positive $%
\star $-functional $l$ such that 
\begin{equation*}
l\left( b\right) =l\left( a\cdot b\right) =l\left( b\cdot c\right) =l\left(
a\cdot b\cdot c\right) =0,\forall a,c\in \mathfrak{a\quad }\Rightarrow \quad
b=0,
\end{equation*}
and it is bounded with respect to $l$ in the sense 
\begin{equation}
\left\| b\right\| =\sup \left\{ \left\| a\cdot b\cdot c\right\|
_{+}^{-}/\left\| a\right\| ^{-}\left\| c\right\| _{+}:a,c\in \mathfrak{a}%
\right\} <\infty \quad \forall b\in \mathfrak{a,}  \label{3.4}
\end{equation}
where $\left\| a\right\| _{+}^{-}=\left| l\left( a\right) \right| ,\quad
\left\| a\right\| ^{-}=l\left( a\cdot a^{\star }\right) ^{1/2},\quad \left\|
a\right\| _{+}=l\left( a^{\star }\cdot a\right) ^{1/2}$, then it is
four-normed in the above sense. The defined by $l$ semi-norms $\left\| \cdot
\right\| _\nu ^\mu $ are obviously separating, satisfy the inequalities (\ref%
{3.3}), and they also satisfy the $\star $-equalities of the following
definition

\begin{definition}
An associative four-normed $\star $-algebra $\mathfrak{a}$ is called
B*-algebra if it is complete in the uniform topology, induced by the
semi-norms $\left( \left\Vert a\right\Vert _{\nu }^{\mu }\right) _{\nu
=+,\bullet }^{\mu =-,\bullet }$, satisfying the following equalities 
\begin{equation}
\left\Vert a\cdot a^{\star }\right\Vert =\left\Vert a\right\Vert \left\Vert
a^{\star }\right\Vert ,\quad \left\Vert a\cdot a^{\star }\right\Vert
_{+}^{-}=\left\Vert a\right\Vert ^{-}\left\Vert a^{\star }\right\Vert
_{+}\quad \forall a\in \mathfrak{a}..  \label{3.5}
\end{equation}%
The It\^{o} B*-algebra is a B*-algebra with self-adjoint annihilator $%
d_{t}=d_{t}^{\star }$, $a\cdot d_{t}=0=d_{t}\cdot a$, $\forall a\in 
\mathfrak{a}$ called the death of $\mathfrak{a}$, and the semi-norms (\ref%
{3.1}) given by a linear positive $\star $-functional $l\left( a^{\star
}\right) =\overline{l\left( a\right) }$, $l\left( a\cdot a^{\star }\right)
\geq 0,\forall a\in \mathfrak{a}$ normalized as $l\left( \text{$d_{t}$}%
\right) =1$.
\end{definition}

Obviously, any C*-algebra can be considered as a B*-algebra in the above
sense with three trivial semi-norms $\left\| a\right\| _{+}^{-}=\left\|
a\right\| ^{-}=\left\| a\right\| _{+}=0,\forall a\in \mathfrak{a}$.
Moreover, as it follows from the inequalities (\ref{3.3}) for $c=1$, every
unital B*-algebra is a C*-algebra, the three nontrivial semi-norms on which
might be given by a state $l$, normalized as $l\left( 1\right) =\left\|
1\right\| _{+}^{-}$. However, if a B*-algebra $\mathfrak{a}$ contains only
approximative identity $e_i\nearrow 1$, and $\left\| e_i\right\|
_{+}^{-}\longrightarrow \infty $, it is a proper dense sub-algebra of its
C*-algebraic completion w.r.t. the norm $\left\| \cdot \right\| $.

Note that the use of the term B*-algebra, the obsolete name for the
C*-algebras, in a more general sense is not contradictive, and will never
make a confusion in the context of It\^{o} algebras, as there is no It\^{o}
algebra which is simultaneously a C*-algebra. Every C*-algebra $\mathcal{A}%
\subseteq \mathcal{B}\left( \mathcal{H}\right) $ with a ciclic vector $\eta
\in \mathcal{H}$ can be embedded into a faithful It\^{o} B*-algebra $%
\mathfrak{a}=\mathbb{C}d_{t}+\mathfrak{x}$ , where the subspace $\mathfrak{x}%
\subseteq \mathfrak{a}$ identified with the factor-algebra $\mathfrak{a}/%
\mathbb{C}d_{t}$, is the operator C*-algebra $\mathcal{A}$. This can be done
by $a=\left( \eta |A\eta \right) d_{t}+x$ with $x=A$, $x^{\star }=A^{\dagger
}$ such that $\mathfrak{x}=\mathcal{A}$, $l\left( x\right) =0$\textrm{,} $%
\left\Vert x\right\Vert _{+}=\left\Vert A\eta \right\Vert $, $\left\Vert
x\right\Vert ^{-}=\left\Vert A^{\dagger }\eta \right\Vert $. However in the
general case the zero mean algebra $\mathfrak{x}=\left\{ a\in \mathfrak{a}%
:l\left( a\right) =0\right\} $ equipped with the factor-product 
\begin{equation}
aa^{\star }=a\cdot a^{\star }-l\left( a\cdot a^{\star }\right) d_{t}=x\cdot
x^{\star }-l\left( x\cdot x^{\star }\right) d_{t}=xx^{\star },  \label{3.7}
\end{equation}%
where $x=a-l\left( a\right) d_{t}$, is not a C*-algebra for an arbitrary It%
\^{o} B*-algebra $\mathfrak{a}$, although it is complete with respect to the
C*-semi-norm $\left\Vert x\right\Vert =\left\Vert a\right\Vert $, jointly
with two Hilbert semi-norms $\left\Vert x\right\Vert _{+}=\left\Vert
a\right\Vert _{+}$, $\left\Vert x\right\Vert ^{-}=\left\Vert a\right\Vert
^{-}$.

Thus in order to classify the It\^{o} B*-algebras we should study the
structure of the zero-mean algebras $\mathfrak{x}=\left\{ x=a-l\left(
a\right) :a\in \mathfrak{a}\right\} $ with the normes, given by the
semi-scalar products $\left\langle x^{\star }|x\right\rangle
_{+}=\left\langle x|x^{\star }\right\rangle ^{-}$, defining the It\^{o}
algebras $\mathfrak{a}=\left\{ \alpha d_{t}+x\right\} $ with $l\left(
a\right) =\alpha $ by 
\begin{equation}
a\cdot a^{\star }=\left\langle x|x\right\rangle ^{-}d_{t}+xx^{\star }=x\cdot
x^{\star },\quad a^{\star }\cdot a=\left\langle x|x\right\rangle
_{+}d_{t}+x^{\star }x=x^{\star }\cdot x.  \label{3.8}
\end{equation}

As was proved in the first section, every It\^{o} B*-algebra is
algebraically and isometrically isomorphic to a closed $\star $-subalgebra $%
\mathfrak{a}\ni D_{t}$ of the simple vacuum B*-algebra $\mathfrak{b}\left( 
\mathcal{H}\right) $ associated with the GNS Hilbert space $\mathcal{H}$ by (%
\ref{0.9}) so that $\boldsymbol{x}=\boldsymbol{i}\left( x\right) $ for each $%
x=a-l\left( a\right) d_{t}$ is given by the tripple 
\begin{equation*}
\hat{x}=i\left( x\right) ,\quad x\rangle _{+}=k\left( x\right) ,\quad
\langle x=k^{\dagger }\left( x\right) ,
\end{equation*}%
with $\left\Vert x\right\Vert =\left\Vert i\left( x\right) \right\Vert $, $%
\left\Vert x\right\Vert _{+}=\left\Vert k\left( x\right) \right\Vert $, $%
\left\Vert x\right\Vert ^{-}=\left\Vert k^{\dagger }\left( x\right)
\right\Vert $. Thus the factor-algebra $\mathfrak{x}=\mathfrak{a}/\mathbb{C}%
d_{t}$ can be identified with a closed $\star $-subalgebra of the B*-algebra 
$\mathcal{H}\oplus \mathcal{A}\oplus \mathcal{H}^{\dagger }$ equipped with
the product 
\begin{equation}
x^{\star }x=\hat{x}^{\dagger }x\rangle _{+}\oplus \hat{x}^{\dagger }\hat{x}%
\oplus \langle x|\hat{x},\quad xx^{\star }=\hat{x}|x\rangle ^{-}\oplus \hat{x%
}\hat{x}^{\dagger }\oplus \langle x\hat{x}^{\dagger },  \label{3.9}
\end{equation}%
for $x=$ $x\rangle _{+}\oplus \hat{x}\oplus \langle x\in \mathcal{H}\oplus 
\mathcal{A}\oplus \mathcal{H}^{\dagger }\ni |x\rangle ^{-}\oplus \hat{x}%
^{\dagger }\oplus \langle x|=x^{\star }$, where 
\begin{equation*}
|x\rangle ^{-}=\langle x^{\dagger }=x^{\star }\rangle _{+},\quad \langle
x|=x\rangle _{+}^{\dagger }=\langle x^{\star }.
\end{equation*}

Let $P$ be the maximal orthoprojector in $\mathcal{H}$ for which $i\left(
a\right) P=0$ for all $a\in \mathfrak{a}$ such that $E=I-P$ is the support
for the operator algebra $\mathcal{A}=i\left( \mathfrak{a}\right) $. We
shall prove that the $\star $-projections 
\begin{equation}
\pi \left( x\right) =Px\rangle _{+}\oplus \langle xP,\quad \varepsilon
\left( x\right) =Ex\rangle _{+}\oplus \hat{x}\oplus \langle xE  \label{3.10}
\end{equation}
define the homomorphisms of $\mathfrak{x}\subseteq \mathcal{H}\oplus 
\mathcal{A}\oplus \mathcal{H}^{\dagger }$ respectively onto the $\star $%
-ideal $\mathfrak{y}=\left\{ y\in \mathfrak{x}:\mathfrak{a}y=0=y\mathfrak{a}%
\right\} $ and the maximal subalgebra $\mathfrak{z}\subseteq \mathfrak{x}$
having the dense part $\mathfrak{aa}=\left\{ ac:a,c\in \mathfrak{a}\right\} $
w.r.t. any of two Hilbert seminorms. An Ito B*-algebra $\mathfrak{a}$ with
the trivial factor-product $\mathfrak{aa}=\left\{ 0\right\} $ and thus $%
\left\| \cdot \right\| =0$ is called the (general) Brownian algebra. In the
opposite case, when $\mathfrak{aa}$ is dense in $\mathfrak{x}$ w.r.t. any of
the Hilbert semi-norms $\left\| \cdot \right\| ^{-},\left\| \cdot \right\|
_{+\,}$, it is called the (general) L\'{e}vy algebra.

\begin{theorem}
Let $\mathfrak{a}$ be a B*-It\^{o} algebra. Then it is an orthogonal sum $%
\mathfrak{b}+\mathfrak{c}$, $\mathfrak{b}\cdot \mathfrak{c}=\left\{
0\right\} $ of a quantum Brownian B*-algebra $\mathfrak{b}$ and a quantum L%
\'{e}vy B*-algebra $\mathfrak{c}$. This decomposition is unique up to the
ideal $\mathfrak{b}\cap \mathfrak{c}=\mathbb{C} d_{t} $.
\end{theorem}

\begin{proof}
We want to find the orthogonal decomposition 
\begin{equation*}
x=y+z,\quad \left\langle y|z\right\rangle _{+}=0=\left\langle
y|z\right\rangle ^{-},\quad y\in \mathfrak{y},z\in \mathfrak{z}
\end{equation*}
for all $x=a-l\left( a\right) $, $a\in \mathfrak{a}$, and to prove that it
is unique. First note that the condition $\mathfrak{a}y=\left\{ 0\right\} =y%
\mathfrak{a} $ for the elements $y\in \mathfrak{y}$ is equivalent to $%
\left\| y\right\| =0$ with the right and left orthogonality of $y\in 
\mathfrak{x}$ to $\mathfrak{aa}\subseteq \mathfrak{x}$: 
\begin{equation*}
\left\langle y|a^{\star }c\right\rangle _{+}=0=\left\| y\right\|
\Leftrightarrow \left\langle ay|c\right\rangle _{+}=0=\left\langle
ya|c\right\rangle _{+}\quad \forall a,c\in \mathfrak{a}
\end{equation*}
\begin{equation*}
\left\langle y|ac^{\star }\right\rangle ^{-}=0=\left\| y\right\|
\Leftrightarrow \left\langle yc|a\right\rangle ^{-}=0=\left\langle
cy|a\right\rangle ^{-}\quad \forall a,c\in \mathfrak{a}
\end{equation*}
So the $\star $-ideal $\mathfrak{y}\subseteq \mathfrak{x}$ with all trivial
products $xy=0=yx$ is defined as the closed $\star $-subspace $\left\{ y\in 
\mathfrak{x}_0:\left\langle y|\mathfrak{aa}\right\rangle _{+}=0=\left\langle
y|\mathfrak{aa}\right\rangle ^{-}\right\} $ of the Hilbert space $\left\{
x\in \mathfrak{x}:\left\| x\right\| =0\right\} $ with the $\star $-invariant
scalar product 
\begin{equation*}
\left\langle z^{\star }|x\right\rangle =\left\langle z^{\star
}|x\right\rangle _{+}+\left\langle x|z^{\star }\right\rangle
^{-}=\left\langle z|x^{\star }\right\rangle ^{-}+\left\langle x^{\star
}|z\right\rangle _{+}=\left\langle x^{\star }|z\right\rangle .
\end{equation*}
Let us prove that the symmetric orthogonal projection $\left\langle y|\pi
\left( x\right) \right\rangle =\left\langle y|x\right\rangle $ of $x\in 
\mathfrak{x}$ onto $\mathfrak{y}\ni y$ is defined as in (\ref{3.10}).
Indeed, 
\begin{eqnarray*}
\left\langle y|x\right\rangle &=&\left\langle y|x\right\rangle
_{+}+\left\langle x|y\right\rangle ^{-}=\left\langle y|Px\right\rangle
_{+}+\left\langle xP|y\right\rangle ^{-} \\
&=&\left\langle y|\pi \left( x\right) \right\rangle _{+}+\left\langle \pi
\left( x\right) |y\right\rangle ^{-}=\left\langle y|\pi \left( x\right)
\right\rangle ,
\end{eqnarray*}
where $P$ is the left orthoprojector $\left\langle y|Px\right\rangle
_{+}=\left\langle y|x\right\rangle _{+}$ in $\mathcal{H}$ onto $\mathcal{G}%
_{+}=\mathfrak{y}\rangle _{+}$, defining the right orthoprojector $P^{\prime
}$ onto $\mathcal{G}^{-}=\mathcal{G}_{+}^{\dagger }$ by $\mathcal{H}%
^{\dagger }\ni \langle x\mapsto \langle xP$. The $\star $-projection $\pi
\left( x^{\star }\right) =\pi \left( x\right) ^{\star }$ is the homorphism
as $\pi \left( xz\right) =0$ due to $\left\langle y|xz\right\rangle =0$ for
all $x,z\in \mathfrak{x}$, $y\in \mathfrak{y}$, and $\varepsilon \left(
x\right) =x-\pi \left( x\right) $, defined in (\ref{3.10}) is also a $\star $%
-homomorphism, $\varepsilon \left( xz\right) =xz$, for all $x,z\in \mathfrak{%
x}$ as its kernel $\varepsilon ^{-1}\left( 0\right) =\left\{ x\in \mathfrak{x%
}:\varepsilon \left( x\right) =0\right\} $ is the $\star $-ideal $\mathfrak{y%
}$. Thus the range $\mathfrak{z}=\varepsilon \left( \mathfrak{x}\right) $ is
a closed $\star $-subalgebra, which is left and right orthogonal to $%
\mathfrak{y}$: 
\begin{equation*}
\left\langle y|\varepsilon \left( x\right) \right\rangle _{+}=\left\langle
y|Ex\right\rangle _{+}=0=\left\langle xE|y\right\rangle ^{-}=\left\langle
\varepsilon \left( x\right) |y\right\rangle ^{-},\quad \forall x\in 
\mathfrak{x},y\in \mathfrak{y}.
\end{equation*}
Due to $\mathcal{E}_{+}=E\mathfrak{x}\rangle _{+}$ contains $\mathfrak{aa}%
\rangle _{+}=\mathcal{A}\mathfrak{x}\rangle _{+}$ as a dense part in the
Hilbert space $\mathcal{H}$, the algebra $\mathfrak{z}$ contains $\mathfrak{%
aa}$ as a left (right) dense part with respect to the left (right) Hilbert
seminorm $\left\| \cdot \right\| _{+}$ ($\left\| \cdot \right\| ^{-}$.)
Obviously $\mathfrak{z}\subseteq \mathfrak{x}$ is uniquely defined as the
maximal such $\star $-subalgebra, and $\varepsilon $ is uniquely defined as
the representation $\mathfrak{z}$ of the quotient algebra $\mathfrak{x}/%
\mathfrak{y}$.

Thus $a=b+c$, $bc=0$ for all $a\in \mathfrak{a}$, where $b=\beta d_{t}+y\in 
\mathfrak{b}$ define the Brownian B*-subalgebra with the fundamental
representation $\mathfrak{b}\subseteq \mathbb{C}\oplus \mathcal{G}_{+}\oplus 
\mathcal{G}^{-}$, where $\mathcal{G}_{+}=Pk\left( \mathfrak{a}\right) $, $%
\mathcal{G}^{-}=k^{\dagger }\left( \mathfrak{a}\right) P$, and $c=\left(
\alpha -\beta \right) d_{t}+z\in \mathfrak{c}$ define the L\'{e}vy
B*-subalgebra, having the fundamental representation $\mathfrak{c}\subseteq 
\mathbb{C\oplus \mathcal{E}}_{+}\oplus \mathcal{E}^{-}\mathbb{\oplus }%
\mathcal{A}$ with non-degenerated operator algebra $\mathcal{A}=i\left( 
\mathfrak{a}\right) $, left and right represented on $\mathcal{E}%
_{+}=Ek\left( \mathfrak{a}\right) $ and $\mathcal{E}^{-}=k^{\dagger }\left( 
\mathfrak{a}\right) E$.
\end{proof}

\begin{example}
The commutative multiplication table $\mathrm{d}w_{i}\mathrm{d}\bar{w}%
_{k}=\delta _{k}^{i}\mathrm{d}t$ for the complex It\^{o} differentials $%
\mathrm{d}w_{k}=\mathrm{d}\bar{w}_{-k}$, $k=\mathbb{Z}$ of the Fourier
amplitudes 
\begin{equation*}
w_{k}\left( t\right) =\int_{-\pi }^{\pi }e^{jk\theta }W\left( t,d\theta
\right) ,
\end{equation*}
\begin{equation*}
\mathrm{d}W\left( t,\Delta \right) \mathrm{d}W\left( t,\Delta ^{\prime
}\right) =\frac{1}{2\pi }\left| \Delta ^{\prime }\cap \Delta \right| \mathrm{%
d}t
\end{equation*}
for an orthogonal Wiener measure $W\left( t,\Delta \right) $ on $\Delta
\subseteq \left[ -\pi ,\pi \right] \ni \theta $ of the normalized intensity
can be generalized in the following way.

Let $\rho _{k}>0$, $k\in \mathbb{Z}$ be a self-inverse family of spectral
eigen-values $\rho _{-k}=\rho _{k}^{-1}$ for a positive-definite
(generalized) periodic function 
\begin{equation*}
\lambda \left( \theta \right) =\sum_{k=-\infty }^{\infty }\rho
_{k}e^{jk\theta },
\end{equation*}%
\begin{equation*}
\left[ \bar{\lambda}\ast \lambda \right] \left( \theta \right) :=\frac{1}{%
2\pi }\int_{-\pi }^{\pi }\bar{\lambda}\left( \theta -\phi \right) \lambda
\left( \phi \right) \mathrm{d}\phi =\delta \left( \theta \right) .
\end{equation*}%
The generalized multiplication table $d_{i}\cdot d_{k}^{\star }=\rho
_{k}\delta _{k}^{i}d_{t}$ for abstract infinitesimals $d_{k}=d_{-k}^{\star }$%
, $k\in \mathbb{Z}$ is obviously non-commutative for all $k$ with $\lambda
_{k}\neq 1$. The $\star $-semigroup $\left\{ 0,\text{$d_{t}$},d_{k}:k\in 
\mathbb{Z}\right\} $ generates an infinite dimensional It\^{o} algebra $%
\left( \mathfrak{b},l\right) $ of a quantum Wiener periodic motion on $\left[
-\pi ,\pi \right] $ as it is the second order nilpotent algebra $\mathfrak{a}
$ of $a=\alpha d_{t}+y$, $y=\sum \eta ^{k}d_{k}$ with $y^{\star }=\sum \eta
_{k}d_{k}^{\star }$ and $l\left( y\right) =0$ for all $\eta ^{k}=\bar{\eta}%
_{k}\in \mathbb{C}$. It is a Brownian B*-algebra with closed involution on
the complex space $\mathcal{D}$ of all $\eta $ given by all complex
sequences $\eta =\left( \eta _{k}\right) _{k\in \mathbb{Z}}$ with 
\begin{equation*}
\left\Vert \eta \right\Vert ^{-}=\left( \sum \left\vert \eta _{k}\right\vert
^{2}\rho _{k}\right) ^{1/2}=\left\Vert \eta ^{\star }\right\Vert _{+}<\infty 
\text{.}
\end{equation*}%
The operator representation of $\mathfrak{a}$ in Fock space is defined by
the forward differentials of $\Lambda \left( t,a\right) =\alpha t\mathrm{I}%
+\sum \eta ^{k}\hat{w}_{k}\left( t\right) $, where 
\begin{equation*}
\hat{w}_{k}\left( t\right) =\hat{v}_{k}^{-}\left( t\right) +\hat{v}%
_{k}^{+}\left( t\right) ,
\end{equation*}%
$\,$%
\begin{equation*}
\hat{v}_{k}^{-}\left( t\right) =\rho _{-k}^{1/2}\int_{-\pi }^{\pi
}e^{jk\theta }\Lambda _{-}\left( t,d\theta \right) ,\quad \hat{v}%
_{k}^{+}\left( t\right) =\rho _{k}^{1/2}\int_{-\pi }^{\pi }e^{jk\theta
}\Lambda ^{+}\left( t,d\theta \right) 
\end{equation*}%
are given by the annihilation and creation measures in Fock space over
square integrable functions on $\mathbb{R}_{+}\times \left[ -\pi ,\pi \right]
$ with the standard multiplication table 
\begin{equation*}
\mathrm{d}\Lambda _{-}\left( t,\Delta \right) \mathrm{d}\Lambda ^{+}\left(
t,\Delta ^{\prime }\right) =\frac{1}{2\pi }\left\vert \Delta ^{\prime }\cap
\Delta \right\vert I\mathrm{d}t,
\end{equation*}%
\begin{equation*}
\mathrm{d}\Lambda ^{+}\left( t,\Delta ^{\prime }\right) \mathrm{d}\Lambda
_{-}\left( t,\Delta \right) =0.
\end{equation*}
\end{example}

\begin{example}
The commutative multiplication table $\mathrm{d}m_i\mathrm{d}\bar{m}%
_k=\delta _k^i\mathrm{d}t+\mathrm{d}m_{i-k}$ for the complex It\^{o}
differentials $\mathrm{d}m_k=\mathrm{d}\bar{m}_{-k}$, $k=\mathbb{Z}$ of the
Fourier amplitudes $m_k\left( t\right) =\int_{-\pi }^\pi e^{jk\theta
}M\left( t,d\theta \right) $, 
\begin{equation*}
\mathrm{d}M\left( t,\Delta ^{\prime }\right) =\frac 1{2\pi }\left| \Delta
^{\prime }\cap \Delta \right| \mathrm{d}t+\mathrm{d}M\left( t,\Delta
^{\prime }\cap \Delta \right) ,\quad \Delta ,\Delta ^{\prime }\subseteq %
\left[ -\pi ,\pi \right]
\end{equation*}
for the standard compensated Poisson measure $M\left( t,\Delta \right) $ can
be generalized in the following way.

Let $G$ be a discrete locally compact group, and $G\ni g\mapsto \lambda
_{g}\in \mathbb{C}$ be a positive-definite summable function, $\lambda
_{g^{-1}}=\bar{\lambda}_{g}$, which is self-inverse in the convolutional
sense $\left[ \bar{\lambda}\ast \lambda \right] _{g}=\sum_{h}\bar{\lambda}%
_{gh^{-1}}\lambda _{h}=\delta _{g}^{1}$. The generalized multiplication
table for abstract infinitesimals 
\begin{equation*}
d_{g}=d_{-g}^{\star },\quad d_{g}\cdot d_{h}^{\star }=\lambda _{gh^{-1}}%
\text{$d_{t}$}+d_{gh^{-1}},\quad g,h\in G
\end{equation*}%
is obviously associative and commutative if $G$ is Abelian, $G\simeq \mathbb{%
Z}$, but it is non-commutative for non Abelian $G$ even if $\lambda
_{k}=\delta _{k}^{1}$ as in the above case. The $\star $-semigroup $\left\{
0,\text{$d_{t}$},d_{g}:g\in G\right\} $ generates an infinite dimensional It%
\^{o} algebra $\left( \mathfrak{a},l\right) $ of a quantum compensated
Poisson motion on the spectrum $\Omega $ of the group $G$ as $\mathfrak{a}$
is the sum of $\mathbb{C}d_{t}$ and the unital group algebra $\mathcal{D}$
of $z=\sum \zeta ^{g}d_{g}$ with involution $z^{\star }=\sum \zeta
_{g}d_{g}^{\star }$ and $l\left( z\right) =0$ for all $\zeta ^{g}=\bar{\zeta}%
_{g}\in \mathbb{C}$. It is a L\'{e}vy B*-algebra with closed involution on
the complex space $\mathcal{D}$ of all complex sequences $\zeta =\left(
\zeta _{g}\right) _{g\in G}$ with 
\begin{equation*}
\left\Vert \zeta \right\Vert ^{-}=\left( \sum \left( \bar{\zeta}\ast \tilde{%
\zeta}\right) _{g}^{2}\lambda _{g}\right) ^{1/2}=\left\Vert \zeta ^{\star
}\right\Vert _{+}<\infty ,
\end{equation*}%
where $\tilde{\zeta}=\left( \zeta _{g^{-1}}\right) _{g\in G}$. The operator
integral representation of $\mathfrak{a}$ is defined in Fock space over the
Hilbert space of square integrable function on $\mathbb{R}_{+}$ with values
in the direct integral $\mathcal{H}=\int_{\Omega }^{\oplus }\mathcal{H}%
\left( \omega \right) \mathrm{d}\pi _{\omega }$ of Hilbert spaces $\mathcal{H%
}\left( \omega \right) $ for the spectral decomposition 
\begin{equation*}
\lambda _{g}=\int_{\Omega }\mathrm{Tr}\left[ \rho _{\omega }U_{g}\left(
\omega \right) \right] \mathrm{d}\pi _{\omega }
\end{equation*}%
w.r.t. the unitary irreducible representations $U_{g}\left( \omega \right) $
of $G$ and the Plansherel measure $\mathrm{d}\pi _{\omega }$. It is given by
the forward differentials of $\Lambda \left( t,a\right) =\alpha t\mathrm{I}%
+\sum \zeta ^{g}\hat{m}_{g}\left( t\right) $, with $\,$%
\begin{eqnarray*}
&&\hat{m}_{g}\left( t\right)  \\
&=&\int_{\Omega }\mathrm{Tr}_{\mathcal{H}\left( \omega \right) }\left\{
U_{g}\left( \omega \right) \left[ \rho _{\omega }^{1/2}\Lambda _{-}\left(
t,d\omega \right) +\Lambda ^{+}\left( t,d\omega \right) \rho _{\omega
}^{1/2}+\Lambda \left( t,d\omega \right) \right] \right\} 
\end{eqnarray*}%
defined by the standard annihilation, creation and exchange operator-valued
measures in Fock space with a measurable family $\left( \rho _{\omega
}\right) _{\omega \in \Omega }$ of positive operators in $\mathcal{H}\left(
\omega \right) $ having the integrable w.r.t. $\mathrm{d}\pi _{\omega }$
traces $\mathrm{Tr}\rho _{\omega }<\infty $.
\end{example}


\begin{thebibliography}{9}
\bibitem{Ito51} It\^{o}, K. On a Formula Concerning Stochastic
Differentials. Nagoya Math . J., \textbf{3}, pp. 55-65, 1951.

\bibitem{Grd91} Gardiner, C.W. Quantum Noise. Springer-Verlag, 1991.

\bibitem{3} Belavkin, V.\thinspace P. Chaotic States and Stochastic
Integration in Quantum Systems. Russian Math.\thinspace Survey, \textbf{47},
(1), pp.\thinspace 47--106, 1992.

\bibitem{4} Belavkin V.P. Kernel Representations of $\star $-semigroups
Associated with Infinitely Divisible States. Quantum Probability and Related
Topics, Vol VII, pp 31--50, 1992.

\bibitem{1} Hudson, R.\thinspace L.\ and Parthasarathy, K.\thinspace R.
Quantum It\^{o}'s Formula and Stochastic Evolution. Comm.\thinspace
Math.\thinspace Phys., \textbf{93}, pp.\thinspace 301--323, 1984.

\bibitem{5} Takesaki, M. J. Functional Analysis \textbf{9}, p 306, 1972.
\end{thebibliography}
\end{document}